\documentclass{amsart}
\usepackage{amssymb, epsf, amscd}
\allowdisplaybreaks
\newtheorem{theorem}[subsubsection]{Theorem}               
\newtheorem{acknowledgements}[subsubsection]{Acknowledgements}

\newtheorem{lemma}[subsubsection]{Lemma}                   
\newtheorem{corollary}[subsubsection]{Corollary}

\theoremstyle{definition}                      
\newtheorem{definition}[subsubsection]{Definition}         
\newtheorem{example}[subsubsection]{Example}                
               
\theoremstyle{remark}
\def\rtb{
\begin{flushright}
\boxed{}
\end{flushright}}
\newtheorem{remark}[subsubsection]{Remark}                 
               
\numberwithin{equation}{subsection}
\def\Arc[#1]{
\ifcase#1
\qbezier[25](0.966,-0.259)(1.04,0)(0.966,0.259)
\or
\qbezier[25](0.966,0.259)(0.897,0.518)(0.707,0.707)
\or
\qbezier[25](0.707,0.707)(0.518,0.897)(0.259,0.966)
\or
\qbezier[25](0.259,0.966)(0,1.04)(-0.259,0.966)
\or
\qbezier[25](-0.259,0.966)(-0.518,0.897)(-0.707,0.707)
\or
\qbezier[25](-0.707,0.707)(-0.897,0.518)(-0.966,0.259)
\or
\qbezier[25](-0.966,0.259)(-1.04,0)(-0.966,-0.259)
\or
\qbezier[25](-0.966,-0.259)(-0.897,-0.518)(-0.707,-0.707)
\or
\qbezier[25](-0.707,-0.707)(-0.518,-0.897)(-0.259,-0.966)
\or
\qbezier[25](-0.259,-0.966)(0,-1.04)(0.259,-0.966)
\or
\qbezier[25](0.259,-0.966)(0.518,-0.897)(0.707,-0.707)
\or
\qbezier[25](0.707,-0.707)(0.897,-0.518)(0.966,-0.259)
\fi}

\newcommand{\fig}[1]
        {\raisebox{-0.5\height}%
                 {\epsfbox{#1}}%
        }

\newcommand\lbb[1]{\label{#1} 
                   }                                    
\def\Zset{\mathbb{Z}}       

\def\CK{\mbox{K}}
\def\CM{{\sf M}}

\def\CF{{\sf F}}
\def\CG{{\sf G}}
\def\be{\beta}

\def\be{\begin{equation}}
\def\fe{\end{equation}}
\def\KI{Z}

\def\CA{A}
\def\Picture#1{
\begin{picture}(2,1)(0,0)
#1
\end{picture}
}

\begin{document}

\title{Branched cyclic covers and finite type invariants}
\author{Andrew Kricker}
\address{Department of Mathematical and Computing Science, Tokyo Institute of Tech
nology}
\email{kricker@is.titech.ac.jp} 
\date{First version: September, 1999,\ \ Current version: February, 2000.}


\begin{abstract}
This work identifies a class of moves on knots which translate to 
$m$-equivalences of the associated
$p$-fold branched cyclic covers, for a fixed $m$ and any $p$
(with respect to the Goussarov-Habiro filtration). 
These moves are applied to give a flexible (if specialised) 
construction of knots
for which the Casson-Walker-Lescop invariant (for example) 
of their 
$p$-fold branched cyclic covers may be readily calculated, for any choice
of $p$.

In the second part of this paper,
these operations are illustrated by some theorems
concerning the relationship of knot invariants obtained from
finite type three-manifold invariants,
via the
branched cyclic covering construction, 
with the finite type theory of knots.
\end{abstract}

\maketitle

\tableofcontents

\def\clasp{
\Picture{
\qbezier(0,1)(0,0.6)(0,0.6)
\qbezier(0,-1)(0,-0.6)(0,0.2)
\qbezier(0,0.4)(0.4,0.4)(0.4,0)
\qbezier(0,0.4)(-0.4,0.4)(-0.4,0)
\qbezier(0.4,0)(0.4,-0.2)(0.2,-0.36)
\qbezier(-0.4,0)(-0.4,-0.2)(-0.2,-0.36)
\qbezier(0.4,0)(0.6,0)(1,0)
}
}

\def\claspf{
\Picture{
\qbezier(0,1)(0,0)(0,-1)
\put(1,0){$\claspc$}
}
}

\def\claspff{
\Picture{
\qbezier[5](0,1)(0,0.5)(0,0)
\qbezier(0,0)(0,-0.5)(0,-1)
}
}

\def\claspffp{
\Picture{
\qbezier(0,0)(0,-0.5)(0,-1)
}
}

\def\claspfft{
\Picture{
\qbezier(0,1)(0,0.5)(0,0)
\qbezier[5](0,0)(0,-0.5)(0,-1)
}
}

\def\claspfftp{
\Picture{
\qbezier(0,1)(0,0.5)(0,0)
}
}

\def\claspc{
\Picture{
\qbezier(0,0.4)(0.4,0.4)(0.4,0)
\qbezier(0,0.4)(-0.4,0.4)(-0.4,0)
\qbezier(0.4,0)(0.4,-0.4)(0,-0.4)
\qbezier(-0.4,0)(-0.4,-0.4)(-0,-0.4)
\qbezier(0.4,0)(0.6,0)(1,0)
}
}

\def\claspr{
\Picture{
\qbezier(-0,-1)(-0,-0.6)(-0,-0.6)
\qbezier(-0,--1)(-0,--0.6)(-0,-0.2)
\qbezier(-0,-0.4)(-0.4,-0.4)(-0.4,-0)
\qbezier(-0,-0.4)(--0.4,-0.4)(--0.4,-0)
\qbezier(-0.4,-0)(-0.4,--0.2)(-0.2,--0.36)
\qbezier(--0.4,-0)(--0.4,--0.2)(--0.2,--0.36)
\qbezier(-0.4,-0)(-0.6,-0)(-1,-0)
}
}

\def\capl{
\Picture{
\qbezier(0,-1)(0,0)(-1,0)
\qbezier(-2,-1)(-2,0)(-1,0)
}}

\def\cupl{
\Picture{
\qbezier(0,1)(0,0)(-1,0)
\qbezier(-2,1)(-2,0)(-1,0)
}}

\def\caplf{
\Picture{
\qbezier[5](-2,-1)(-2,-1.5)(-2,-2)
\qbezier[10](0,-1)(0,0)(-1,0)
\qbezier[10](-2,-1)(-2,0)(-1,0)
}}

\def\cuplf{
\Picture{
\qbezier[5](-2,1)(-2,1.5)(-2,2)
\qbezier[10](0,1)(0,0)(-1,0)
\qbezier[10](-2,1)(-2,0)(-1,0)
}}

\def\leftbox{
\Picture{
\put(1,0){$\clasp$}
\put(1,-2){$\clasp$}
\put(1,2){$\capl$}
\put(1,-4){$\cupl$}
\put(-1,-3){\line(0,1){4}}
}}

\def\leftboxn{
\Picture{
\put(1,0){$\claspff$}
\put(1,-2){$\claspfft$}
\put(1,2){$\caplf$}
\put(1,-4){$\cuplf$}
\put(-1,-2){\line(0,1){2}}
}}

\def\leftboxnb{
\Picture{
\put(1,0){$\claspff$}
\put(1,-2){$\claspfftp$}
\put(1,2){$\caplf$}
\put(-1,-2){\line(0,1){2}}
\qbezier(-1,-2)(-1,-3)(0,-3)
\qbezier(1,-2)(1,-3)(0,-3)
}}

\def\leftboxnc{
\Picture{
\put(1,0){$\claspffp$}
\put(1,-2){$\claspfft$}
\put(1,-4){$\cuplf$}
\put(-1,-2){\line(0,1){2}}
\qbezier(-1,0)(-1,1)(0,1)
\qbezier(1,0)(1,1)(0,1)
}}

\def\reltripa{
\Picture{
\put(-9,5){\circle*{0.7}}
\put(-4,4){$\leftboxn$}
\qbezier(-5,3)(-7,3)(-9,5)
\put(-9,5){\line(0,1){2}}
\put(-9,5){\line(-1,0){2}}
\put(-6,-2){$(M,L_A)$}}
}

\def\reltripb{
\Picture{
\put(-9,5){\circle*{0.7}}
\put(-4,4){$\leftboxnb$}
\qbezier(-5,3)(-7,3)(-9,5)
\put(-9,5){\line(0,1){2}}
\put(-9,5){\line(-1,0){2}}
\put(-6,-2){$(M,L_B)$}
}
}

\def\reltripc{
\Picture{
\put(-9,5){\circle*{0.7}}
\put(-4,4){$\leftboxnc$}
\qbezier(-5,3)(-7,3)(-9,5)
\put(-9,5){\line(0,1){2}}
\put(-9,5){\line(-1,0){2}}
\put(-6,-2){$(M,L_C)$}
}
}

\def\leftboxfa{
\Picture{
\put(1,0){$\claspf$}
\put(1,-2){$\clasp$}
\put(1,2){$\capl$}
\put(1,-4){$\cupl$}
\put(-1,-3){\line(0,1){4}}
}}

\def\leftboxfb{
\Picture{
\put(1,0){$\clasp$}
\put(1,-2){$\claspf$}
\put(1,2){$\capl$}
\put(1,-4){$\cupl$}
\put(-1,-3){\line(0,1){4}}
}}

\def\leftboxf{
\Picture{
\put(1,0){$\claspc$}
\put(1,-2){$\claspc$}
}}

\def\leftboxp{
\Picture{
\put(1,0){$\clasp$}
\put(1,-2){$\clasp$}
\put(1,2){$\capl$}
\put(1,-4){$\cupl$}
\qbezier(-1,-3)(-3,-3)(-5,-1)
}}


\def\capr{
\Picture{
\qbezier(-0,--1)(-0,-0)(--1,-0)
\qbezier(--2,--1)(--2,-0)(--1,-0)
}}

\def\cupr{
\Picture{
\qbezier(-0,-1)(-0,-0)(--1,-0)
\qbezier(--2,-1)(--2,-0)(--1,-0)
}}

\def\rightbox{
\Picture{
\put(-1,-0){$\claspr$}
\put(-1,--2){$\claspr$}
\put(-1,-2){$\capr$}
\put(-1,--4){$\cupr$}
\put(1,3){\line(-0,-1){4}}
}}

\def\vertclasp{
\Picture{
\qbezier(-1,0)(-0.8,0)(-0.6,0)
\qbezier(--1,0)(--0.6,-0)(-0.2,-0)
\qbezier(-0.4,0)(-0.4,0.4)(-0,0.4)
\qbezier(-0.4,0)(-0.4,-0.4)(--0,-0.4)
\qbezier(-0,0.4)(--0.2,0.4)(--0.36,0.2)
\qbezier(-0,-0.4)(--0.2,-0.4)(--0.36,-0.2)
\qbezier(-0,0.4)(-0,0.6)(-0,1)
}}

\def\vertclaspb{
\Picture{
\qbezier(-1,0)(-0.8,0)(-0.6,0)
\qbezier(--1,0)(--0.6,-0)(-0.2,-0)
\qbezier(-0.4,0)(-0.4,-1.4)(--0,-1.4)
\qbezier(-0,-1.4)(--0.2,-1.4)(--0.36,-0.2)
\qbezier(-0.4,0)(-0.4,3)(-1,3)
\qbezier(0.4,0.2)(0.4,3)(1,3)
}}

\def\vertclaspud{
\Picture{
\qbezier(--1,0)(--0.8,0)(--0.6,0)
\qbezier(---1,0)(---0.6,-0)(--0.2,-0)
\qbezier(--0.4,0)(--0.4,-0.4)(--0,-0.4)
\qbezier(--0.4,0)(--0.4,--0.4)(---0,--0.4)
\qbezier(--0,-0.4)(---0.2,-0.4)(---0.36,-0.2)
\qbezier(--0,--0.4)(---0.2,--0.4)(---0.36,--0.2)
\qbezier(--0,-0.4)(--0,-0.6)(--0,-1)
}}

\def\ygraph{
\Picture{
\put(0,-3.75){$\vertclaspf$}
\put(-3.45,6){$\vertclaspfu$}
\put(3.45,6){$\vertclaspfu$}
\put(0,0){\circle{1.5}}
\qbezier(0.75,0)(3.75,1)(3.75,3)
\qbezier(-0.75,0)(-3.75,1)(-3.75,3)
\qbezier(-0.45,0.5)(-3.25,1.5)(-3.15,3)
\qbezier(0.45,0.5)(3.25,1.5)(3.15,3)
}}

\def\vertclaspf{
\Picture{
\qbezier(-0.3,1)(-0.3,2)(-0.3,3)
\qbezier(0.3,1)(0.3,2)(0.3,3)
\qbezier(-2,-0.2)(-3,-0.2)(-3,-4)
\qbezier(-2,0.2)(-3.6,0.2)(-3.6,-4)
\qbezier(2,-0.2)(3,-0.2)(3,-4)
\qbezier(2,0.2)(3.6,0.2)(3.6,-4)
\qbezier(1,-0.2)(1,-1)(0,-1)
\qbezier(-1,0)(-1,1)(0,1)
\qbezier(1,0.2)(1,1)(0,1)
\qbezier(-1,0)(-1,-1)(0,-1)
\qbezier(-2,-0.2)(-1,-0.2)(-1,-0.2)
\qbezier(--2,-0.2)(--0.6,-0.2)(-0.36,-0.2)
\qbezier(-2,0.2)(-1,0.2)(-1,0.2)
\qbezier(--2,0.2)(--0.6,0.2)(-0.36,0.2)
\qbezier(-0.4,0)(-0.4,0.4)(-0,0.4)
\qbezier(-0.4,0)(-0.4,-0.4)(--0,-0.4)
\qbezier(-0,0.4)(--0.2,0.4)(--0.36,0.2)
\qbezier(-0,-0.4)(--0.2,-0.4)(--0.36,-0.2)
}}

\def\vertclaspfu{
\Picture{
\qbezier(-0.3,-1)(-0.3,-2)(-0.3,-3)
\qbezier(0.3,-1)(0.3,-2)(0.3,-3)
\qbezier(-2,--0.2)(-3,--0.2)(-3,--4)
\qbezier(-2,-0.2)(-3.6,-0.2)(-3.6,--4)
\qbezier(2,--0.2)(3,--0.2)(3,--4)
\qbezier(2,-0.2)(3.6,-0.2)(3.6,--4)
\qbezier(1,--0.2)(1,--1)(0,--1)
\qbezier(-1,-0)(-1,-1)(0,-1)
\qbezier(1,-0.2)(1,-1)(0,-1)
\qbezier(-1,-0)(-1,--1)(0,--1)
\qbezier(-2,--0.2)(-1,--0.2)(-1,--0.2)
\qbezier(--2,--0.2)(--0.6,--0.2)(-0.36,--0.2)
\qbezier(-2,-0.2)(-1,-0.2)(-1,-0.2)
\qbezier(--2,-0.2)(--0.6,-0.2)(-0.36,-0.2)
\qbezier(-0.4,-0)(-0.4,-0.4)(-0,-0.4)
\qbezier(-0.4,-0)(-0.4,--0.4)(--0,--0.4)
\qbezier(-0,-0.4)(--0.2,-0.4)(--0.36,-0.2)
\qbezier(-0,--0.4)(--0.2,--0.4)(--0.36,--0.2)
}}

\def\extvertclasp{
\Picture{
\put(0,0){$\vertclasp$}
\qbezier(-4,0)(-1,0)(-1,0)
\qbezier(4,0)(1,0)(1,0)
}}

\def\verthook{
\Picture{
\qbezier(-1,0)(-1,-2)(0,-2)
\qbezier(0,-2)(1,-2)(1,-0.2)
\qbezier(-4,0)(-1.25,0)(-1.25,0)
\qbezier(4,0)(1.25,0)(-0.75,0)
\qbezier(-1,0)(-1,1)(-1,2)
\qbezier(1,0.2)(1,1)(1,2)
}}

\def\verthookc{
\Picture{
\qbezier(-1,0)(-1,-2)(0,-2)
\qbezier(0,-2)(1,-2)(1,-0.7)
\qbezier(-4,0)(-1.25,0)(-1.25,0)
\qbezier(-4,-0.5)(-1.25,-0.5)(-1.25,-0.5)
\qbezier(4,0)(1.25,0)(-0.75,0)
\qbezier(4,-0.5)(1.25,-0.5)(-0.75,-0.5)
\qbezier(-1,0)(-1,1)(-1,2)
\qbezier(1,0.2)(1,1)(1,2)
}}

\def\verthookp{
\Picture{
\qbezier(-1,0)(-1,-2)(0,-2)
\qbezier(-4,0)(-1.25,0)(-1.25,0)
\qbezier(5,0)(1.25,0)(-0.75,0)
\qbezier(-1,0)(-1,1)(-1,2)
\qbezier(1,0.2)(1,1)(1,2)
}}

\def\halftwist{
\Picture{
\put(0,0){\circle{1.4}}
\put(0,0.7){\line(0,1){0.3}}
\put(0,-0.7){\line(0,-1){0.3}}
\qbezier(0.49,0.49)(0,0)(-0.49,-0.49)
}
}

\def\halftwistp{
\Picture{
\put(0,0){\circle{1.4}}
\put(0.7,0){\line(1,0){0.3}}
\put(-0.7,0){\line(-1,0){0.3}}
\qbezier(0.49,-0.49)(0,0)(-0.49,0.49)
}
}

\def\shade{
\Picture{
\qbezier[20](0,-3)(-2,-5)(-4,-7)
\qbezier[20](1,-3)(-1,-5)(-3,-7)
\qbezier[20](2,-3)(0,-5)(-2,-7)
\qbezier[20](-2,-3)(-4,-5)(-6,-7)
\qbezier[20](-1,-3)(-3,-5)(-5,-7)
\qbezier[20](-0.2,-3)(-2.2,-5)(-4.2,-7)
\qbezier[20](-0.4,-3)(-2.4,-5)(-4.4,-7)
\qbezier[20](-0.6,-3)(-2.6,-5)(-4.6,-7)
\qbezier[20](-0.8,-3)(-2.8,-5)(-4.8,-7)
\qbezier[20](0.2,-3)(-1.8,-5)(-3.8,-7)
\qbezier[20](0.4,-3)(-1.6,-5)(-3.6,-7)
\qbezier[20](0.6,-3)(-1.4,-5)(-3.4,-7)
\qbezier[20](0.8,-3)(-1.2,-5)(-3.2,-7)
\qbezier[20](-2.2,-3)(-4.2,-5)(-6.2,-7)
\qbezier[20](-2.4,-3)(-4.4,-5)(-6.4,-7)
\qbezier[20](-2.6,-3)(-4.6,-5)(-6.6,-7)
\qbezier[20](-2.8,-3)(-4.8,-5)(-6.8,-7)
\qbezier[20](2.2,-3)(0.2,-5)(-1.8,-7)
\qbezier[20](2.4,-3)(0.4,-5)(-1.6,-7)
\qbezier[20](2.6,-3)(0.6,-5)(-1.4,-7)
\qbezier[20](2.8,-3)(0.8,-5)(-1.2,-7)
\qbezier[20](-1.2,-3)(-3.2,-5)(-5.2,-7)
\qbezier[20](-1.4,-3)(-3.4,-5)(-5.4,-7)
\qbezier[20](-1.6,-3)(-3.6,-5)(-5.6,-7)
\qbezier[20](-1.8,-3)(-3.8,-5)(-5.8,-7)
\qbezier[20](1.2,-3)(-0.8,-5)(-2.8,-7)
\qbezier[20](1.4,-3)(-0.6,-5)(-2.6,-7)
\qbezier[20](1.6,-3)(-0.4,-5)(-2.4,-7)
\qbezier[20](1.8,-3)(-0.2,-5)(-2.2,-7)
}}

\def\shadet{
\Picture{
\qbezier[10](1,-3)(-1,-5)(-3,-7)
\qbezier[10](-1,-3)(-3,-5)(-5,-7)
\qbezier[10](-0.2,-3)(-2.2,-5)(-4.2,-7)
\qbezier[10](-0.6,-3)(-2.6,-5)(-4.6,-7)
\qbezier[10](0.2,-3)(-1.8,-5)(-3.8,-7)
\qbezier[10](0.6,-3)(-1.4,-5)(-3.4,-7)
\qbezier[10](-2.2,-3)(-4.2,-5)(-6.2,-7)
\qbezier[10](-2.6,-3)(-4.6,-5)(-6.6,-7)
\qbezier[10](2.2,-3)(0.2,-5)(-1.8,-7)
\qbezier[10](2.6,-3)(0.6,-5)(-1.4,-7)
\qbezier[10](-1.4,-3)(-3.4,-5)(-5.4,-7)
\qbezier[10](-1.8,-3)(-3.8,-5)(-5.8,-7)
\qbezier[10](1.4,-3)(-0.6,-5)(-2.6,-7)
\qbezier[10](1.8,-3)(-0.2,-5)(-2.2,-7)
}}

\def\leg{
\Picture{
\put(0,-3){$\extvertclasp$}
\qbezier(0,-2)(0,-1)(0,-1)
\put(0,-1){\circle*{0.7}}
\qbezier(0,-1)(2.5,3)(4,3)
\qbezier(0,-1)(-2.5,3)(-4,3)
\put(0,0){$\shade$}
}}

\def\legc{
\Picture{
\qbezier(0,-2)(0,-1)(0,-1)
\put(0,-1){\circle*{0.7}}
\qbezier(0,-1)(2,1)(3,2)
\qbezier(0,-1)(-2,1)(-3,2)
\put(0,-3){$\halftwist$}
\qbezier(0,-4)(0,-6)(0,-6)
\put(-2,-8){$(M,L_A)$}
}}

\def\legcc{
\Picture{
\qbezier(0,-2)(0,-1)(0,-1)
\put(0,-1){\circle*{0.7}}
\qbezier(0,-1)(2,1)(3,2)
\qbezier(0,-1)(-2,1)(-3,2)
\qbezier(0,-2)(0,-6)(0,-6)
\put(-2,-8){$(M,L_B)$}
}}

\def\legccp{
\Picture{
\qbezier(0,-2)(0,-1)(0,-1)
\put(0,-1){\circle*{0.7}}
\qbezier(0,-1)(2,1)(4,3)
\qbezier(0,-1)(-2,1)(-4,3)
\qbezier(0,-2)(0,-6)(0,-9)
}}

\def\legp{
\Picture{
\put(0,-3){$\extvertclasp$}
\qbezier(0,-2)(0,-1)(0,-1)
\put(0,-1){\circle*{0.7}}
\qbezier(0,-1)(2,1)(3,2)
\qbezier(0,-1)(-2,1)(-3,2)
\put(-3,2){\circle*{0.7}}
\qbezier(-3,2)(-3,3)(-3,5)
\put(-3,5){\circle*{0.7}}
\qbezier(-3,5)(-2.5,5.5)(-2,6)
\qbezier(-3,5)(-3.5,5.5)(-4,6)
\put(-6,2){\circle*{0.7}}
\qbezier(-6,2)(-6.5,2.5)(-7,3)
\qbezier(-6,2)(-6.5,1.5)(-7,1)
\qbezier(-3,2)(-4,2)(-6,2)
\put(-1.5,-6.5){$(K,L)$}
}}

\def\legppb{
\Picture{
\put(0,-3){$\extvertclasp$}
\put(4,-4){$K$}
\qbezier(0,-2)(0,-1)(0,-1)
\put(0,-1){\circle*{0.7}}
\qbezier(0,-1)(2,1)(3,2)
\qbezier(0,-1)(-2,1)(-3,2)
}}

\def\leglink{
\Picture{
\put(0,-3){$\extvertclasp$}
\put(4,-4){$K$}
\qbezier(0,-2)(0,-1)(0,2)
}}

\def\polya{
\Picture{
\put(-2.8,2.2){$x_1$}
\put(-4.3,-5){\vector(0,1){4}}
\qbezier(-7,-3)(-5,-3)(-5,-3)
\put(0,0){$\shade$}
\put(0,-3){$\extvertclasp$}
\put(5,-2){$K$}
\put(-2,1){\vector(1,-1){0.01}}
\put(2,1){\vector(1,1){0.01}}
\qbezier(0,-2)(0,-1)(0,-1)
\put(0,-1){\circle*{0.7}}
\qbezier(0,-1)(2,1)(3,2)
\qbezier(0,-1)(-2,1)(-3,2)
\put(-2.8,-10){$1-t^{-x^i}$}
}}

\def\polyb{
\Picture{
\put(-2.8,2.2){$x_1$}
\put(-4.3,-5){\vector(0,1){4}}
\qbezier(-7,-3)(-5,-3)(-5,-3)
\put(0,0){$\shade$}
\put(0,-3){$\extvertclasp$}
\put(5,-2){$K$}
\put(-2,1){\vector(-1,1){0.01}}
\put(2,1){\vector(-1,-1){0.01}}
\qbezier(0,-2)(0,-1)(0,-1)
\put(0,-1){\circle*{0.7}}
\qbezier(0,-1)(2,1)(3,2)
\qbezier(0,-1)(-2,1)(-3,2)
\put(-2.8,-10){$1-t^{x^i}$}
}}

\def\addlegs{
\Picture{
\put(-3,-3){
\vector(1,0){0.01}}
\put(0,-3){$\extvertclasp$}
\put(5,-2){$K$}
\put(-2,1){\vector(-1,1){0.01}}
\put(2,1){\vector(-1,-1){0.01}}
\qbezier(0,-2)(0,-1)(0,-1)
\put(0,-1){\circle*{0.7}}
\qbezier(0,-1)(2,1)(3,2)
\qbezier(0,-1)(-2,1)(-3,2)
}}

\def\legpp{
\Picture{
\put(-3,5){\circle*{0.7}}
\qbezier(-3,5)(-2.5,5.5)(-2,6)
\qbezier(-3,5)(-3.5,5.5)(-4,6)
\put(-6,2){\circle*{0.7}}
\qbezier(-6,2)(-6.5,2.5)(-7,3)
\qbezier(-6,2)(-6.5,1.5)(-7,1)
\put(-3,2){\circle*{0.7}}
\qbezier(-3,2)(-3,3)(-3,5)
\qbezier(-3,2)(-4,2)(-6,2)
\qbezier(-3,2)(0,-1)(3,2)
\qbezier(-4,-3)(0,-3)(4,-3)
\put(-1.5,-6.5){$(K,L')$}
}}

\def\legppc{
\Picture{
\put(-3,2){\circle*{0.7}}
\qbezier(-3,2)(-3,3)(-3,5)
\qbezier(-3,2)(-4,2)(-6,2)
\put(3,2){\circle*{0.7}}
\qbezier(3,2)(3,3)(3,5)
\qbezier(3,2)(4,2)(6,2)
\qbezier(-3,2)(-2,1)(0,1)
\qbezier(3,2)(2,1)(0,1)
\put(-1.5,-6.5){$(M,L_A)$}
}}

\def\yot{
\Picture{
\put(0,0){$\cross$}
\put(-1.6,-0.75){$\crossr$}
\put(-1,1.4){\line(0,1){3}}
\put(-0.6,-1.15){\line(0,-1){3}}
}}

\def\crossh{
\Picture{
\qbezier(1,-1.4)(0.5,-1.4)(0,-1.4)
\qbezier(1,-1.4)(1.4,-1.4)(1.4,-1)
\qbezier(1,-0.6)(1.4,-0.6)(1.4,-1)
\qbezier(0,-1.4)(-0.4,-1.4)(-0.4,-1)
\qbezier(-0.4,-1)(-0.4,-0.6)(0,-0.6)
}}

\def\crosshb{
\Picture{
\qbezier(-1,--1.4)(-0.5,--1.4)(-0,--1.4)
\qbezier(-1,--1.4)(-1.4,--1.4)(-1.4,--1)
\qbezier(-1,--0.6)(-1.4,--0.6)(-1.4,--1)
\qbezier(-0,--1.4)(--0.4,--1.4)(--0.4,--1)
\qbezier(--0.4,--1)(--0.4,--0.6)(-0,--0.6)
}}

\def\break{
\Picture{
\qbezier(-1.4,0)(-3,0)(-3,0)
\qbezier(3,0.15)(1.2,0.15)(1.2,0.15)
\put(-0.85,0.95){$\crossh$}
\put(0.8,-0.75){$\crosshb$}
\put(-6,2){\circle*{0.7}}
\qbezier(-6,2)(-6,3)(-6,5)
\qbezier(-6,2)(-7,2)(-9,2)
\put(6,2){\circle*{0.7}}
\qbezier(6,2)(6,3)(6,5)
\qbezier(6,2)(7,2)(9,2)
\qbezier(-6,2)(-5,0)(-3,0)
\qbezier(6,2)(5,0.15)(3,0.15)
}}

\def\breakb{
\Picture{
\qbezier(-3,0)(-3,0)(3,0)
\put(-6,2){\circle*{0.7}}
\qbezier(-6,2)(-6,3)(-6,5)
\qbezier(-6,2)(-7,2)(-9,2)
\put(6,2){\circle*{0.7}}
\qbezier(6,2)(6,3)(6,5)
\qbezier(6,2)(7,2)(9,2)
\qbezier(-6,2)(-5,0)(-3,0)
\qbezier(6,2)(5,0)(3,0)
}}

\def\brek{
\Picture{
\qbezier(-1.4,0)(-3,0)(-3,0)
\qbezier(3,0.15)(1.2,0.15)(1.2,0.15)
\put(-0.85,0.95){$\crossh$}
\put(0.8,-0.75){$\crosshb$}
}}

\def\brekb{
\Picture{
\put(-3,0){$\brek$}
\put(3,0.15){$\brek$}
}}

\def\legpppb{
\Picture{
\put(4,-4){$K$}
\qbezier(-3,2)(0,-1)(3,2)
\qbezier(-4,-3)(0,-3)(4,-3)
}}

\def\legppppb{
\Picture{
\qbezier(-3,2)(0,-1)(3,2)
}}

\def\legppp{
\Picture{
\put(-3,5){\circle*{0.7}}
\qbezier(-3,5)(-2.5,5.5)(-2,6)
\qbezier(-3,5)(-3.5,5.5)(-4,6)
\put(-6,2){\circle*{0.7}}
\qbezier(-6,2)(-6.5,2.5)(-7,3)
\qbezier(-6,2)(-6.5,1.5)(-7,1)
\put(-3,2){\circle*{0.7}}
\qbezier(-3,2)(-3,3)(-3,5)
\qbezier(-3,2)(-4,2)(-6,2)
\qbezier(-3,2)(-1,2)(-1,-1)
\put(0,-3){$\verthook$}
\qbezier(1,-1)(1,1)(3,2)
\put(-1.5,-6.5){$(K,L'')$}
}}

\def\legpppc{
\Picture{
\put(-3,2){\circle*{0.7}}
\qbezier(-3,2)(-3,3)(-3,5)
\qbezier(-3,2)(-4,2)(-6,2)
\put(3,2){\circle*{0.7}}
\qbezier(3,2)(3,3)(3,5)
\qbezier(3,2)(4,2)(6,2)
\qbezier(-3,2)(-1,2)(-1,-1)
\put(0,-3){$\verthookc$}
\qbezier(1,-1)(1,1)(3,2)
\put(-1.5,-6.5){$(M,L_B)$}
}}

\def\legged{
\Picture{
\put(-4,4){$\leftbox$}
\put(4,2){$\rightbox$}
\qbezier(-5,3)(-7,3)(-9,3)
\qbezier(5,3)(7,3)(9,3)
\qbezier(-2,4)(0,4)(2,4)
\put(0,-4){$\verthookp$}
\put(1,-1){$\shade$}
\qbezier(-2,2)(-1,2)(-1,0)
\qbezier(2,2)(1,2)(1,0)
\put(-1,-1){$\halftwist$}
\qbezier(1,0)(1,-1)(1,-2)
}}

\def\leggedy{
\Picture{
\put(-4,4){$\leftbox$}
\put(4,2){$\rightbox$}
\put(1,2){$\vertclaspud$}
\put(1,1){\line(0,-1){6}}
\qbezier(-5,3)(-7,3)(-9,7)
\qbezier(5,3)(7,3)(9,7)
\qbezier(-2,4)(0,4)(2,4)
\put(-1,2){$\halftwistp$}
}}

\def\leggedp{
\Picture{
\put(-4,4){$\leftbox$}
\put(4,2){$\rightbox$}
\qbezier(-5,3)(-7,3)(-9,5)
\qbezier(5,3)(7,3)(9,3)
\qbezier(-2,4)(0,4)(2,4)
\put(6,-8.75){$\check{F}$}
\qbezier(-2,2)(-1,2)(-1,0)
\qbezier(2,2)(1,2)(1,0)
\put(-1,-1){$\halftwist$}
\put(-0.75,1.75){B}
\qbezier(1,0)(1,-1)(1,-2)
\put(0,-5){$\modhookp$}
\qbezier(3,-7)(3,-5.5)(2,-4.5)
\qbezier(1,-2)(1,-3.5)(2,-4.5)
\put(3,-4.5){C}
\put(-9,5){\circle*{0.7}}
\put(-8,2){Y}
\put(-0.4,5){A}
\put(-9,5){\line(0,1){2}}
\put(-9,5){\line(-1,0){2}}
}}

\def\tripa{
\Picture{
\put(-9,2){Y$^{a_1}$}
\put(2,3.75){A$^{a_1}$}
\put(2,1.75){B$^{a_1}$}
\put(-9,5){\circle*{0.7}}
\put(-4,4){$\leftbox$}
\qbezier(-5,3)(-7,3)(-9,5)
\qbezier(-2,4)(0,4)(1,4)
\qbezier(-2,2)(-1,2)(1,2)
\put(-9,5){\line(0,1){2}}
\put(-9,5){\line(-1,0){2}}
\put(-6,-2){$\widetilde{L}^{(a_1,\ldots,a_s)}$}}
}

\def\tripb{
\Picture{
\put(-9,2){Y$^{a_1}$}
\put(2,3.75){A$^{a_1}$}
\put(2,1.75){B$^{a_1}$}
\put(-9,5){\circle*{0.7}}
\put(-4,4){$\leftboxfa$}
\qbezier(-5,3)(-7,3)(-9,5)
\qbezier(-1,4)(0,4)(1,4)
\qbezier(-2,2)(-1,2)(1,2)
\put(-9,5){\line(0,1){2}}
\put(-9,5){\line(-1,0){2}}
\put(-6,-2){${R}^{(a_1,\ldots,a_s)}$}
}}

\def\tripc{
\Picture{
\put(-9,2){Y$^{a_1}$}
\put(2,3.75){A$^{a_1}$}
\put(2,1.75){B$^{a_1}$}
\put(-9,5){\circle*{0.7}}
\put(-4,4){$\leftboxfb$}
\qbezier(-5,3)(-7,3)(-9,5)
\qbezier(-2,4)(0,4)(1,4)
\qbezier(-1,2)(-1,2)(1,2)
\put(-9,5){\line(0,1){2}}
\put(-9,5){\line(-1,0){2}}
\put(-6,-2){${S}^{(a_1,\ldots,a_s)}$}
}}

\def\leggedpp{
\Picture{
\put(-4,4){$\leftboxf$}
\put(4,2){$\rightbox$}
\qbezier(5,3)(7,3)(9,3)
\qbezier(-2,4)(0,4)(2,4)
\qbezier(-2,2)(-1,2)(-1,0)
\qbezier(2,2)(1,2)(1,0)
\put(-1,-1){$\halftwist$}
\qbezier(1,0)(1,-1)(1,-2)
\put(0,-5){$\modhookp$}
\qbezier(3,-7)(3,-5.5)(2,-4.5)
\qbezier(1,-2)(1,-3.5)(2,-4.5)
}}

\def\pierce{
\Picture{
\qbezier(-1,5)(-1,-2)(-1,-2)
\qbezier[28](-1,-2)(-1,-4)(-1,-6)
\put(0.5,3){$\shadet$}
\qbezier(-5,0)(-1.2,0)(-1.2,0)
\qbezier(-0.8,0)(5,0)(5,0)
}}

\def\modhook{
\Picture{
\qbezier(-1.4,1)(-1.4,0.5)(-1.4,0)
\qbezier(-1.4,1)(-1.4,1.4)(-1,1.4)
\qbezier(-0.6,1)(-0.6,1.4)(-1,1.4)
\qbezier(-1.4,0)(-1.4,-0.4)(-1,-0.4)
\qbezier(-1,-0.4)(-0.6,-0.4)(-0.6,0)
\qbezier(-1.2,1)(-0.6,1)(-0.4,1)
\qbezier(-1.2,0)(-0.6,0)(-0.4,0)
\qbezier(-1.6,1)(-1.8,1)(-1.8,1)
\qbezier(-1.6,0)(-1.8,0)(-1.8,0)
\qbezier(-0.4,0)(0.2,0)(0.2,-2)
\qbezier(-0.4,1)(1.2,1)(1.2,-2)
\qbezier(-1.8,0)(-2.4,0)(-2.4,-2)
\qbezier(-1.8,1)(-3.4,1)(-3.4,-2)
\qbezier(0.2,-2)(0.2,-3)(-0.4,-3)
\qbezier(-3.4,-2)(-3.4,-3)(-3.8,-3)
\qbezier(-2.4,-2)(-2.4,-3)(-2,-3)
\qbezier(1.2,-2)(1.2,-3)(1.6,-3)
\qbezier(-5,-2)(-5,-2)(-3.4,-2)
\qbezier(-2.4,-2)(0,-2)(0.2,-2)
\qbezier(1.2,-2)(2.8,-2)(2.8,-2)
\qbezier(-1,1.4)(-1,3)(-1,3)
\qbezier[20](-1.6,0.5)(-2.9,0.5)(-2.9,-2)
\qbezier[20](-0.4,0.5)(0.7,0.5)(0.7,-2)
\qbezier[6](-1.2,0.5)(-0.8,0.5)(-0.4,0.5)
\qbezier[20](0.7,-2)(0.7,-4)(-1,-4)
\qbezier[20](0.7,-2)(0.7,-4)(-1,-4)
\qbezier[20](-2.9,-2)(-2.9,-4)(-1,-4)
\qbezier[28](-1,-4)(-1,-6.25)(-1,-8.5)
}}

\def\modhookp{
\Picture{
\qbezier(-1.4,1)(-1.4,0.5)(-1.4,0)
\qbezier(-1.4,1)(-1.4,1.4)(-1,1.4)
\qbezier(-0.6,1)(-0.6,1.4)(-1,1.4)
\qbezier(-1.4,0)(-1.4,-0.4)(-1,-0.4)
\qbezier(-1,-0.4)(-0.6,-0.4)(-0.6,0)
\qbezier(-1.2,1)(-0.6,1)(-0.4,1)
\qbezier(-1.2,0)(-0.6,0)(-0.4,0)
\qbezier(-1.6,1)(-1.8,1)(-1.8,1)
\qbezier(-1.6,0)(-1.8,0)(-1.8,0)
\qbezier(-0.4,0)(0.2,0)(0.2,-2)
\qbezier(-0.4,1)(1.2,1)(1.2,-2)
\qbezier(-1.8,0)(-2.4,0)(-2.4,-2)
\qbezier(-1.8,1)(-3.4,1)(-3.4,-2)
\qbezier(0.2,-2)(0.2,-3)(-0.4,-3)
\qbezier(-3.4,-2)(-3.4,-3)(-3.8,-3)
\qbezier(-2.4,-2)(-2.4,-3)(-2,-3)
\qbezier(1.2,-2)(1.2,-3)(1.6,-3)
\qbezier(-9,-2)(-5,-2)(-3.4,-2)
\qbezier(-2.4,-2)(0,-2)(0.2,-2)
\qbezier(1.2,-2)(2.8,-2)(7,-2)
\qbezier(-1,1.4)(-1,3)(-1,3)
\qbezier[20](-1.6,0.5)(-2.9,0.5)(-2.9,-2)
\qbezier[20](-0.4,0.5)(0.7,0.5)(0.7,-2)
\qbezier[6](-1.2,0.5)(-0.8,0.5)(-0.4,0.5)
\qbezier[20](0.7,-2)(0.7,-4)(-1,-4)
\qbezier[20](0.7,-2)(0.7,-4)(-1,-4)
\qbezier[20](-2.9,-2)(-2.9,-4)(-1,-4)
\qbezier[28](-1,-4)(-1,-6)(-1,-6)
\qbezier[20](-1,-6)(-1,-8)(1,-8)
\qbezier[20](3,-6)(3,-8)(1,-8)
\qbezier[32](3,-6)(3,-4)(3,-2)
}}

\def\cross{
\Picture{
\qbezier(-1.4,1)(-1.4,0.5)(-1.4,0)
\qbezier(-1.4,1)(-1.4,1.4)(-1,1.4)
\qbezier(-0.6,1)(-0.6,1.4)(-1,1.4)
\qbezier(-1.4,0)(-1.4,-0.4)(-1,-0.4)
\qbezier(-1,-0.4)(-0.6,-0.4)(-0.6,0)
}}

\def\crossr{
\Picture{
\qbezier(--1.4,1)(--1.4,0.5)(--1.4,0)
\qbezier(--1.4,1)(--1.4,1.4)(--1,1.4)
\qbezier(--0.6,1)(--0.6,1.4)(--1,1.4)
\qbezier(--1.4,0)(--1.4,-0.4)(--1,-0.4)
\qbezier(--1,-0.4)(--0.6,-0.4)(--0.6,0)
}}

\def\seifcross{
\Picture{
\put(0,0){$\cross$}
\put(-1.6,-0.75){$\crossr$}
\put(-1,1.4){\line(0,1){3}}
\put(-0.6,-1.15){\line(0,-1){3}}
\put(-0.5,2.5){$B^{i+1}$}
\put(-0.2,-3){$C^{i}$}
}}

\def\seifcrossb{
\Picture{
\put(0,0){$\cross$}
\put(-1.6,-0.75){$\crossr$}
\put(-1,1.4){\line(0,1){3}}
\put(-0.6,-1.15){\line(0,-1){3}}
\put(-0.5,2.5){Y$_a$}
\put(-0.2,-3){Y$_b$}
}}

\def\ycomp{
\Picture{
\put(5,4.8){\circle{2}}
\put(-5,4.8){\circle{2}}
\put(5,4.8){\circle{1}}
\put(-5,4.8){\circle{1}}
\put(0,0){\circle{1.2}}
\qbezier(-0.3,-0.55)(-0.3,-1)(-0.3,-4)
\qbezier(0.3,-0.55)(0.3,-1)(0.3,-4)
\qbezier(0.58,-0.1)(2.6,1.9)(4.6,3.9)
\qbezier(-0.59,-0.1)(-2.6,1.9)(-4.6,3.9)
\qbezier(0.4,0.4)(2.4,2.4)(4.2,4.2)
\qbezier(-0.4,0.4)(-2.4,2.4)(-4.2,4.2)
\put(0,-5){\circle{2}}
\put(0,-5){\circle{1}}
\put(0,-5){\circle{1}}
}}

\def\actycompa{
\Picture{
\put(5,4.8){\circle{2}}
\put(-5,4.8){\circle{2}}
\put(5,4.8){\circle{1}}
\put(-6,4.6){\vector(0,-1){0.01}}
\put(-5.5,5){\vector(0,1){0.01}}
\put(-1,-5.2){\vector(0,-1){0.01}}
\put(-0.5,-4.8){\vector(0,1){0.01}}
\put(4,4.6){\vector(0,-1){0.01}}
\put(4.5,5){\vector(0,1){0.01}}
\put(-5,4.8){\circle{1}}
\put(0,0){\circle{1.2}}
\put(-0.3,-2){\vector(0,-1){0.01}}
\put(0.3,-2){\vector(0,1){0.01}}
\qbezier(-0.3,-0.55)(-0.3,-1)(-0.3,-4)
\qbezier(0.3,-0.55)(0.3,-1)(0.3,-4)
\qbezier(0.58,-0.1)(2.6,1.9)(4.6,3.9)
\put(2.6,1.9){\vector(1,1){0.01}}
\put(-2.6,1.9){\vector(1,-1){0.01}}
\qbezier(-0.59,-0.1)(-2.6,1.9)(-4.6,3.9)
\qbezier(0.4,0.4)(2.4,2.4)(4.2,4.2)
\put(2.4,2.4){\vector(-1,-1){0.01}}
\put(-2.4,2.4){\vector(-1,1){0.01}}
\qbezier(-0.4,0.4)(-2.4,2.4)(-4.2,4.2)
\put(0,-5){\circle{2}}
\put(0,-5){\circle{1}}
\put(0,-5){\circle{1}}
}}

\def\actycompb{
\Picture{
\put(5,4.8){\circle{2}}
\put(-5,4.8){\circle{2}}
\put(5,4.8){\circle{1}}
\put(-5,4.8){\circle{1}}
\put(0,0){\circle{1.2}}
\qbezier(-0.3,-0.55)(-0.3,-1)(-0.3,-4)
\qbezier(0.3,-0.55)(0.3,-1)(0.3,-4)
\qbezier(0.58,-0.1)(2.6,1.9)(4.6,3.9)
\qbezier(-0.59,-0.1)(-2.6,1.9)(-4.6,3.9)
\qbezier(0.4,0.4)(2.4,2.4)(4.2,4.2)
\qbezier(-0.4,0.4)(-2.4,2.4)(-4.2,4.2)
\put(0,-5){\circle{2}}
\put(0,-5){\circle{1}}
\put(0,-5){\circle{1}}
}}

\def\actycomp{
\Picture{
\put(5,4.8){\circle{2}}
\put(-5,4.8){\circle{2}}
\put(0,0){\circle*{1.2}}
\qbezier(0,-0.6)(0,-4)(0,-4)
\qbezier(0,0)(2.25,2)(4.5,4)
\qbezier(0,0)(-2.25,2)(-4.5,4)
\put(0,-5){\circle{2}}
}}

\def\chord{
\Picture{
\put(-1.5,2){$K$}
\put(4.5,2){$K$}
\qbezier(0,1)(0,2)(0,3)
\qbezier(0,-1)(0,-2)(0,-3)
\qbezier(4,1)(4,2)(4,3)
\qbezier(4,-1)(4,-2)(4,-3)
\put(0,0){$\clasp$}
\qbezier(1,0)(2,0)(3,0)
\put(4,0){$\claspr$}}
}

\def\strut{
\Picture{
\put(-3,0){$\vertclasp$}
\put(3,0){$\vertclasp$}
\qbezier[32](-6,0)(0,0)(6,0)
\put(4,-2){$K$}
\put(0,6){\circle*{1}}
\qbezier(-3,1)(-3,4)(0,6)
\qbezier(3,1)(3,4)(0,6)
\qbezier(0,6)(0,7)(0,9)
}}

\def\thetaA{
\Picture{
\put(-3,0){$\vertclasp$}
\qbezier(-3,1)(-3,2)(-3,3)
\put(3,0){$\vertclasp$}
\qbezier(3,1)(3,2)(3,3)
\qbezier[20](-3,0)(0,0)(3,0)
\put(-3,3){\circle*{1}}
\put(0,5){\circle*{1}}
\put(0,7){\circle*{1}}
\qbezier(0,5)(0,6)(0,7)
\put(3,3){\circle*{1}}
\qbezier[5](-2,1.5)(0,1.5)(2,1.5)
\qbezier(-5,3)(-3,3)(-2,3)
\qbezier(5,3)(3,3)(2,3)
\qbezier[16](-2,3)(0,3)(2,3)
\qbezier(-5,3)(-7,3)(-7,5)
\qbezier(5,3)(7,3)(7,5)
\put(-7,5){\circle*{1}}
\put(7,5){\circle*{1}}
\qbezier(-7,5)(0,5)(7,5)
\qbezier(-5,7)(-7,7)(-7,5)
\qbezier(5,7)(7,7)(7,5)
\qbezier(-5,7)(0,7)(5,7)
\qbezier(-4,0)(-6,0)(-6,-2)
\qbezier(4,0)(6,0)(6,-2)
\qbezier(-6,-2)(-6,-4)(0,-4)
\qbezier(6,-2)(6,-4)(0,-4)
\put(-3.5,-1){$\underbrace{\ \ \ \ \ \ \ \ \ \ \ \ \ \ \ \ \ \ \ }_n$}
\put(7.5,6){$\kappa_n$}
\put(6.75,-2){$U$}
}}

\def\thetaTESTp{
\Picture{
\qbezier(-1,3)(-0.5,3)(0,3)
\qbezier(-2,0)(-0.5,0)(0,0)
\qbezier(-4,3)(-3,3)(-2,3)
\qbezier(-4,0)(-3,0)(-2,0)
\qbezier(4,3)(4,3)(5,3)
\put(3,0){$\vertclaspb$}
\put(1,0){$\vertclaspb$}
\qbezier(-1,3)(-1.5,3)(-2,3)
\put(-1,3){\circle*{1}}
\put(-1,7){\circle*{1}}
\qbezier(-1,3)(-1,5)(-1,7)
\qbezier(-5,3)(-3,3)(-2,3)
\qbezier(-5,3)(-7,3)(-7,5)
\qbezier(5,3)(7,3)(7,5)
\qbezier(-5,7)(-7,7)(-7,5)
\qbezier(5,7)(7,7)(7,5)
\qbezier(-5,7)(0,7)(5,7)
\qbezier(-4,0)(-6,0)(-6,-2)
\qbezier(4,0)(6,0)(6,-2)
\qbezier(-6,-2)(-6,-4)(0,-4)
\qbezier(6,-2)(6,-4)(0,-4)
}}

\def\thetaTEST{
\Picture{
\put(-5,5){$a$}
\put(1,5){$b$}
\qbezier[12](-4,5)(-4,4)(-3,4)
\qbezier[12](-2,5)(-2,4)(-3,4)
\qbezier[12](-2,5)(-2,6)(-3,6)
\qbezier[12](2,5)(2,4)(3,4)
\qbezier[12](4,5)(4,4)(3,4)
\qbezier[12](4,5)(4,6)(3,6)
\put(-3,6){\vector(-1,0){0.01}}
\put(3,6){\vector(-1,0){0.01}}
\put(-3.3,-1.5){$\epsilon_1$}
\put(2.7,-1.5){$\epsilon_2$}
\qbezier(-1,3)(-0.5,3)(0,3)
\qbezier(-2,0)(-0.5,0)(0,0)
\put(-3,0){$\vertclasp$}
\qbezier(-3,1)(-3,2)(-3,3)
\put(-3,3){\circle*{1}}
\put(3,0){$\vertclasp$}
\put(1,0){$\vertclaspb$}
\qbezier(3,1)(3,2)(3,3)
\put(3,3){\circle*{1}}
\qbezier(-1,3)(-1.5,3)(-2,3)
\put(-1,3){\circle*{1}}
\put(-1,7){\circle*{1}}
\qbezier(-1,3)(-1,5)(-1,7)
\qbezier(-5,3)(-3,3)(-2,3)
\qbezier(5,3)(3,3)(2,3)
\qbezier(-5,3)(-7,3)(-7,5)
\qbezier(5,3)(7,3)(7,5)
\qbezier(-5,7)(-7,7)(-7,5)
\qbezier(5,7)(7,7)(7,5)
\qbezier(-5,7)(0,7)(5,7)
\qbezier(-4,0)(-6,0)(-6,-2)
\qbezier(4,0)(6,0)(6,-2)
\qbezier(-6,-2)(-6,-4)(0,-4)
\qbezier(6,-2)(6,-4)(0,-4)
}}

\def\thetaTESTq{
\Picture{
\qbezier(-1,3)(-0.5,3)(0,3)
\qbezier(-2,0)(-0.5,0)(0,0)
\put(-3,0){$\vertclasp$}
\qbezier(-3,1)(-3,2)(-3,3)
\put(-3,3){\circle*{1}}
\put(3,0){$\vertclasp$}
\put(1,0){$\vertclaspb$}
\qbezier(3,1)(3,2)(3,3)
\put(3,3){\circle*{1}}
\qbezier(-1,3)(-1.5,3)(-2,3)
\put(-1,3){\circle*{1}}
\put(-1,7){\circle*{1}}
\qbezier(-1,3)(-1,5)(-1,7)
\qbezier(-5,3)(-3,3)(-2,3)
\qbezier(5,3)(3,3)(2,3)
\qbezier(-5,3)(-7,3)(-7,5)
\qbezier(5,3)(7,3)(7,5)
\qbezier(-5,7)(-7,7)(-7,5)
\qbezier(5,7)(7,7)(7,5)
\qbezier(-5,7)(0,7)(5,7)
\qbezier(-4,0)(-6,0)(-6,-2)
\qbezier(4,0)(6,0)(6,-2)
\qbezier(-6,-2)(-6,-4)(0,-4)
\qbezier(6,-2)(6,-4)(0,-4)
}}

\def\thetaTESTr{
\Picture{
\qbezier(-1,3)(-0.5,3)(0,3)
\qbezier(-2,0)(-0.5,0)(0,0)
\put(-3,0){$\vertclasp$}
\qbezier(-3,1)(-3,2)(-3,3)
\put(-3,3){\circle*{1}}
\put(3,0){$\vertclasp$}
\put(1,0){$\vertclaspb$}
\qbezier(3,1)(3,2)(3,3)
\put(3,3){\circle*{1}}
\qbezier(-1,3)(-1.5,3)(-2,3)
\put(-1,3){\circle*{1}}
\put(-1,7){\circle*{1}}
\qbezier(-1,3)(-1,5)(-1,7)
\qbezier(-5,3)(-3,3)(-2,3)
\qbezier(5,3)(3,3)(2,3)
\qbezier(-5,3)(-7,3)(-7,5)
\qbezier(5,3)(7,3)(7,5)
\qbezier(-5,7)(-7,7)(-7,5)
\qbezier(5,7)(7,7)(7,5)
\qbezier(-5,7)(0,7)(5,7)
\qbezier(-4,0)(-6,0)(-6,-2)
\qbezier(4,0)(6,0)(6,-2)
\qbezier(-6,-2)(-6,-4)(0,-4)
\qbezier(6,-2)(6,-4)(0,-4)
}}

\def\thetaTESTs{
\Picture{
\qbezier(-1,3)(-0.5,3)(0,3)
\qbezier(-2,0)(-0.5,0)(0,0)
\put(-3,0){$\vertclasp$}
\qbezier(-3,1)(-3,2)(-3,3)
\put(-3,3){\circle*{1}}
\put(3,0){$\vertclasp$}
\put(1,0){$\vertclaspb$}
\qbezier(3,1)(3,2)(3,3)
\put(3,3){\circle*{1}}
\qbezier(-1,3)(-1.5,3)(-2,3)
\put(-1,3){\circle*{1}}
\put(-1,7){\circle*{1}}
\qbezier(-1,3)(-1,5)(-1,7)
\qbezier(-5,3)(-3,3)(-2,3)
\qbezier(5,3)(3,3)(2,3)
\qbezier(-5,3)(-7,3)(-7,5)
\qbezier(5,3)(7,3)(7,5)
\qbezier(-5,7)(-7,7)(-7,5)
\qbezier(5,7)(7,7)(7,5)
\qbezier(-5,7)(0,7)(5,7)
\qbezier(-4,0)(-6,0)(-6,-2)
\qbezier(4,0)(6,0)(6,-2)
\qbezier(-6,-2)(-6,-4)(0,-4)
\qbezier(6,-2)(6,-4)(0,-4)
}}

\def\thetab{
\Picture{
\put(0,5){\circle*{1}}
\put(0,7){\circle*{1}}
\qbezier(0,5)(0,6)(0,7)
\put(-3,5){\circle*{1}}
\put(3,5){\circle*{1}}
\qbezier(-3,5)(0,5)(3,5)
\qbezier(-3,5)(-3,7)(0,7)
\qbezier(3,5)(3,7)(0,7)
\qbezier(-3,5)(-3,3)(0,3)
\qbezier(3,5)(3,3)(0,3)
\put(3.8,5){$\lambda_{\kappa_n}$}
}}

\def\thetac{
\Picture{
\put(0,5){\circle*{0.3}}
\put(0,7){\circle*{0.3}}
\qbezier[12](0,5)(0,6)(0,7)
\put(-3,5){\circle*{0.3}}
\put(3,5){\circle*{0.3}}
\qbezier[20](-3,5)(0,5)(3,5)
\qbezier[12](-3,5)(-3,7)(0,7)
\qbezier[12](3,5)(3,7)(0,7)
\qbezier[18](-3,5)(-3,3)(0,3)
\qbezier[18](3,5)(3,3)(0,3)
\put(3.8,5){$\kappa$}
}
}

\pagestyle{myheadings}
\markboth{ANDREW KRICKER}{BRANCHED CYCLIC COVERS AND FINITE TYPE INVARIANTS}

\newpage

\section{Introduction}\lbb{sintro}

Let $\CK$ be the free abelian group generated by (oriented) equivalence
classes of pairs $(M,K)$, where $M$ is an integral homology three-sphere
and $K$ is an oriented
knot in it. Let $\CM$ be the free abelian group generated
by homeomorphism classes of three-manifolds.
Let the mapping
\[
\Sigma^p : \CK \rightarrow \CM
\]
be the linear extension of the operation defined on a pair $(M,K)$
as the $p$-fold branched cyclic covering of $M$ branched over $K$.
We are chiefly interested in knots in the three-sphere, for which
the corresponding space will be denoted $\CK(S^3)$, when required.

The subject of this investigation is the set of knot invariants obtained
by composing the branched cyclic cover construction, for a fixed degree,
with finite-type
invariants of three-manifolds (of various descriptions). 

This paper is in two parts. The main theorem of the first part
is Theorem \ref{maintha}, and its enhancement, Theorem
\ref{mainthb}. This theorem describes the effect of surgery on
{\bf complete graph Y-links} on the $p$-fold branched cyclic covers
of a knot, with respect to the Goussarov-Habiro filtration.
A calculus is developed in Lemmas \ref{calclemma} and \ref{calclemmab}
which describes this action. 

In the second part of the paper, we illustrate 
this calculus by exploring
some very general questions about the invariants obtained by composing
projections of the LMO invariant with the operation of taking a $p$-fold
branched cyclic cover. The main theorem of this part is the realisation
theorem, Theorem \ref{notft}, which has the following, possibly unsurprising
corollary:
\begin{corollary}\label{thecorol}
Take a positive integer $p$.
If $v$ is a rational valued three-manifold invariant 
factoring non-trivially through the LMO invariant on integral homology 
spheres, then $v\circ \Sigma^p$ is not a 
finite-type invariant of knots.
\end{corollary} 

It is interesting that we can use
a finite-type property (in one sense) 
to show
that an invariant is not finite-type (in some different sense).

The generality of the above theorem does not seem to follow from the
more descriptive investigations of specific examples and
projections of the LMO invariant that have been recorded:
Hoste gave a formula for the Casson invariant of the
$p$-fold branched cyclic cover over an untwisted double of a knot \cite{Hos}; 
Davidow gave formulae for the case of iterated torus knots \cite{Dav1}, 
and for some $1$-twisted doubles \cite{Dav2};
Mullins calculated the composition of the Casson-Walker invariant
with the 2-fold branched cyclic cover in terms of the Jones polynomial,  when 
the left hand side was well-defined \cite{Mul1,Mul2}; Ishibe used this work
to give formulae for the general case of $m$-twisted doubles \cite{Ish};  
Garoufalidis showed that Mullins formulae was valid in general, using the
Casson-Walker-Lescop invariant \cite{Gar}; and also in that paper 
Garoufalidis initiated the investigation into the LMO invariant on 
branched cyclic covers, describing the form of the result for covers over 
twisted doubles. 

The main motivation for this work was to see whether the techniques we
employed in \cite{K} (an application of clasper theory) had 
wider application. 
The other main 
motivation was to explore a possible approach to
Lev Rozansky's program for a ``finite type theory of knots' complements''
\cite{Roz}. 
We expect the chief interest in this article to be the analysis of a
suggestive family of operations which will be important in the development
of
a more structural relationship between the loop expansion of the Kontsevich
integral, and the degree expansion of $\KI^{LMO}\circ \Sigma^p$. (That is
the informal conjecture we are seeking to motivate with this paper.)

A review of some of the theory of Goussarov-Habiro filtration
is included as an appendix.

\begin{acknowledgements}

The author is supported by a Japan Society for the Promotion of Science
Postdoctoral Fellowship, and thanks Tomotada Ohtsuki and
the Department of Mathematical and Computing Sciences at the 
Tokyo Institute of Technology for their support.
The author would also like to thank
Stavros Garoufalidis, Kazuo Habiro, Hitoshi Murakami 
and Theodore Stanford (and many others) for interesting comments.
\end{acknowledgements}

\section{Part 1: Complete clasper moves and branched cyclic covers}
\label{formsect}

The Goussarov-Habiro filtration \cite{G3,Hab}
is a descending filtration of $\CM$
(see Definition \ref{GHD}):
\[
\CM \supset \CF^Y_1\CM \supset \CF^Y_2\CM \supset \ldots
\]

On the other hand we have clasper calculus, introduced by Habiro \cite{HT,Hab},
certain aspects of which were independently developed by Goussarov \cite{G2}.
This generalised Matveev's ``Borromeo'' move \cite{Mat} 
and Murakami-Nakanishi's ``delta-unknotting'' operation \cite{MN}. 

In this work, we will identify a subclass of clasper operations which the 
branched cyclic covering construction will translate to 
$n$-equivalences of three-manifolds. 
For technical reasons, 
this work will use the language 
of Y-links;
note that our language is slightly different to that
adopted in \cite{GL2}. 

In Appendix \ref{GHapp} 
we will fix the term {\it graph Y-link} to mean
a Y-link decoration of a knot $K$ such that each leaf of every
Y-component links either another leaf, in a ball (see Appendix \ref{GHapp}
for diagrammatic conventions):

\

\vspace{1.5cm}

\hspace{1.5cm}\seifcrossb\hspace{1cm}
or links the knot thus:
\hspace{2.0cm}\leglink

\vspace{2cm}

Note that we also consider graph Y-links which 
may be {\it without legs}. Note that surgery on a graph
Y-link without legs will change the homeomorphism class of the ambient
manifold; surgery on a connected graph Y-link with at least one leg will not.

To every graph Y-link there is canonically associated
(up to sign) a uni-trivalent diagram 
(which here refers to the familiar uni-trivalent
graphs, which may have univalent vertices located on a skeleton; when
we take linear combinations, we introduce
STU, IHX and AS relations). This association is made
by replacing Y-components by trivalent vertices,
replacing leaves linking the knot by legs, and by
joining tines of trivalent vertices by an edge when the corresponding 
leaves were linked. A graph Y-link is called {\it connected} when the 
dashed graph of this associated diagram is connected. The {\it degree} of a 
graph Y-link is half the number
of vertices of the dashed graph of the associated uni-trivalent diagram.

To introduce the moves in question, some definitions are required. 
\begin{definition}
A graph Y-link decorating a knot $K$ is {\it complete} if it is
connected and if every
leg of the associated uni-trivalent diagram leads to a seperate trivalent vertex.
\end{definition}

In particular, we are ruling out the chord:\vspace{1cm}

\hspace{4.5cm}\chord

\vspace{1.25cm}

\begin{flushleft}
and graph Y-links which end in a fork:
\end{flushleft}

\

\vspace{3cm}

\hspace{5.25cm}
\strut\vspace{1cm}

It is important to point out that this is quite a strong restriction: for 
example, a move of such a kind on a link will not affect its Milnor
invariants.

Now we need an alternative measure on graph Y-links.
\begin{definition}
The {\it surplus} of a uni-trivalent diagram is defined to be the number of
trivalent vertices minus the number of univalent vertices of its
dashed graph. The {\it surplus}
of a graph Y-link is defined to be the surplus of its associated uni-trivalent
diagram.
\end{definition}

For example, the surplus of a complete graph Y-link is simply the
grade of the underlying uni-trivalent diagram, after all the legs have been removed.
The main theorem is expressed in these terms.
\begin{theorem}[Main theorem]\label{maintha}

If a knot $K^L$ is obtained from a knot $K$ by surgery on a complete
graph Y-link $L$ of surplus $s$, where $s\geq 2$, then
\[
\Sigma^p_{K^L} - \Sigma^p_K \in \CF^Y_s\CM.
\]
\end{theorem}

There is another important invariant for which a similar property holds.
The following theorem \cite{KGr}, 
and the theory of which it is an application, may
appear in a future manuscript. 

\begin{theorem}\label{chuckin}
If a knot $K^L$ is obtained from a knot $K$ by surgery on a complete
graph Y-link $L$ of surplus $s$, where $s\geq 2$, then
\[
Z^{LMO}(K^L) - Z^{LMO}(K) = \{ \mbox{series of diagrams of surplus $\geq s$} \},
\]
where $Z^{LMO}$ is the Kontsevich-LMO (Le-Murakami-Ohtsuki) invariant
of knots in integral homology three-spheres.
\end{theorem}
The point being that it is true for every grade. We leave the reader to
make their own conjectures.

Returning to the main theorem of this paper, 
the next natural question is how these moves operate on the 
leading term in the space of associated graded quotients. 
Specifically, what is
\[
\ \ \Sigma^p_{K^L} - \Sigma^p_K \in \CG^Y_s\CM \ ?
\]

Lemma \ref{calclemma} and Lemma \ref{calclemmab} calculate this term. 
To report this calculation, let the ambient integral-homology
three-sphere be presented by some surgery link, and isotope $K$ and
the decorating graph Y-link in that three-manifold to appear in the 
complement of that surgery link in the three-sphere.
\begin{definition}
Let $\lambda_L$ be a graph Y-link in $S^3$ obtained by sawing off 
the legs of $L$\vspace{0.5cm}

\[
\hspace{0.5cm}
\legppb\hspace{1.6cm}\longrightarrow\hspace{2.2cm}
\legpppb
\]\vspace{1cm}
\begin{flushleft}
and then forgetting $K$ and the surgery link. 
\end{flushleft}

Let $\lambda_L'$ be the corresponding graph Y-link in $\Sigma^p_K$ induced by the connect-sum of that
three-sphere into $\Sigma^p_K$. Let $(\Sigma^p_K)^{\lambda_L'}$ denote 
the result of surgery on that graph Y-link.
\end{definition}

Lemma \ref{welldefnlem} asserts that $(\Sigma^p_K)^{\lambda_L'}
- \Sigma^p_K$ is well-defined in $\CG^Y_s\CM$. 
Lemma  \ref{calclemma} and Lemma \ref{calclemmab} show that:

\begin{theorem}[Main theorem, part 2]\label{mainthb}
\[
\Sigma^p_{K^L} - \Sigma^p_{K} \propto (\Sigma^p_K)^{\lambda_L'}
- \Sigma^p_K \in \CG^Y_s\CM.
\]
\end{theorem}

\begin{remark}
The definition of $(\Sigma^p_K)^{\lambda_L'}- \Sigma^p_K$ is perhaps
more easily understood in terms of the surjective map from the space of 
marked uni-trivalent graphs of degree $s$ to $\CG^Y_s\CM([M])$. 
In this sense, the first non-vanishing term is proportional to
the image under that map 
of the diagram obtained by removing 
all the legs from a diagram associated
to the initial decoration. The procedure we employ in this paper 
is used so 
that we can be precise, in a natural way, 
about the sign of the element referred to
(see Section \ref{signsect}).
\end{remark}

The following two lemmas can be used to immediately calculate this
term. 
Calculation proceeds in two steps.
The first lemma is used to reduce the number of legs on $K$, so that 
eventually $\Sigma^p_{K^L} - \Sigma^p_K$ is expressed as a linear combination 
of 
differences $\Sigma^p_{K_i} - \Sigma^p_K$ where each $K_i$ has been 
obtained from
$K$ by surgery on a connected graph Y-link without legs. The second
lemma calculates an expression for the terms in that sum. 

\begin{remark}
This procedure is the reason it is more natural in this context
to work in the generality of integral homology three-spheres.
It is interesting that the calculation necessitates this extension, not
unlike the results of \cite{GL2}.
\end{remark}

\begin{lemma}\label{calclemma}

Let $K$ be a knot.
Let $L$, $L'$ and $L''$  denote three
decorations of $K$ by complete graph Y-links
of surplus $s$ that differ in a ball in the following way:\vspace{1.5cm}

\[
\hspace{2cm}
\legp\hspace{1cm},
\]

\vspace{4cm}

\[
\hspace{2.75cm}
\legpp\hspace{2cm}\mbox{and}\hspace{4cm}\legppp
\hspace{1cm}.\ \ \ \ \ 
\]\vspace{2cm}

Then
\[
(\Sigma^p_{K^L} - \Sigma^p_{K}) = (\Sigma^p_{K^{L'}} - \Sigma^p_K)\ -\
(\Sigma^p_{K^{L''}} - \Sigma^p_K)\ \in \CG^Y_s\CM.
\]
\end{lemma}

\

The statement of the next lemma 
requires an embedding of a graph in $M-K$ 
that will be associated
to a decoration of $K$ by a graph Y-link $L$ in the complement of 
$K$. This embedding is constructed by contracting the
clasper graph presenting $L$ onto a spine: call this embedded graph $E_L$.
In the statement of the lemma, a loop on $L$ will mean precisely a
path in $M-K$ that corresponds to a loop on $E_L$.

\begin{lemma}\label{calclemmab}

Let a knot $K^L$
be obtained from a knot $K$ by surgery on some connected 
graph Y-link $L$ without legs of surplus $s$. Then
\[
\Sigma^p_{K^L} - \Sigma^p_K = \left\{
\begin{array}{lll}
p ((\Sigma^p_{K})^{\lambda_L'} - \Sigma^p_{K}) & \mbox{if} & \mbox{every loop on $L$ has zero linking} \\
& & \mbox{number with $K$ modulo $p$}, \\
& & \\
0 & & \mbox{otherwise}.
\end{array}
\right.\ \in\ \CG^Y_s\CM.
\]
\end{lemma}

Intuitively speaking, this is what one 
would expect. The implementation of that intuition in this setting, however, 
is not immediate (Section \ref{calcbsect}).

Let us illustrate this calculus with some examples based on the
most familiar Goussarov-Habiro finite type invariant, the 
Casson-Walker-Lescop invariant, denoted $\lambda_{CWL}$.
The connection will be made with the following lemma.
Extend the symbol $|H_1(M)|$ to be zero in the case that $M$ has
non-vanishing first Betti number.
\begin{lemma}

\

\begin{enumerate}
\item{ $\lambda_{CWL}( \CG^Y_3\CM ) = 0.$ }
\item{ Let some three-manifold $M$ have some graph Y-link $\Theta$ 
embedded in it, whose underlying diagram is the ``theta'' graph. 
Then
\[
\lambda_{CWL}( M^\Theta ) = \lambda_{CWL}(M) \pm 2|H_1(M)|,
\]
where the sign can be determined in a given situation.
}
\end{enumerate}
\end{lemma}
Perhaps the easiest way to see this is to observe the identification
of the co-efficient of the theta graph in the LMO invariant 
with half the Casson-Walker-Lescop
invariant, \cite{LMMO} (see also \cite{GH,HB,L}), and then observe that
the difference of the LMO invariant on a manifold and on that manifold 
surgered along a $\Theta$ is precisely $\pm |H_1(M)| \theta$ plus higher
order terms, where the sign can be determined in a given situation
(see \cite{LeGr} for such first non-vanishing term calculations).

\begin{example}

Fix an integer $p$ and take a knot $K$. Consider some graph $Y$-link 
decoration $\Theta$ of $K$ with underlying diagram the ``theta'' graph, 
and let $K^\Theta$ be the knot (in some integral homology three-sphere) 
obtained from surgery on $\Theta$. Let $\theta$ denote the spatial graph 
in $S^3-K$ associated to this decoration. Let $l_1$
denote a choice of a knot in $S^3-K$ obtained by forgetting 
some edge of $\theta$, and let $l_2$ denote a knot obtained by forgetting
some other choice of edge. Then:

\[
\lambda_{CWL}(\Sigma^p_{K^{\Theta}}) = \lambda_{CWL}(\Sigma^p_K) \pm \left\{
\begin{array}{ll}
2p|H_1(\Sigma^p_K)| & \mbox{if lk$(l_i,K)=0\ $mod p, $i=1,2$\ ,} \\
0 & \mbox{otherwise,}
\end{array}\right.
\]
where the sign can be determined in a given situation. This follows
from Lemma \ref{calclemmab}. 
\rtb

\end{example}

How might we organise such a calculation in generality? That is, where
we have a decoration of a knot by a theta graph with some legs attached.
Such an example follows.
\begin{example}\label{compexamp}
We indicate a basis $\{a,b\}$ 
for the first homology of the associated trivalent graph, and some labels
$\{\epsilon_1,\epsilon_2\}$ have been affixed to the legs. Let $U$ indicate
the unknot, and let the decorating graph Y-link be denoted by $\kappa$.

\

\vspace{2cm}

\[
\thetaTEST\hspace{2.7cm} - \hspace{1.3cm}
U
\]
\vspace{1cm}

Theorem \ref{maintha} indicates that under $\Sigma^p$ this difference
lies in $\CG^Y_s\CM$. To calculate the term we proceed by applying
Lemma \ref{calclemma} to each of the legs. This produces four terms,
which will be recorded by setting each of the $\epsilon_i$ to either
0 or 1. For example the term that corresponds to $(\epsilon_1,\epsilon_2)
=(0,1)$ will be:

\vspace{2.4cm}

\[
\thetaTESTp\hspace{2.7cm} - \hspace{1.3cm} U
\]
\vspace{1.6cm}

Now, according to Lemma \ref{calclemmab}, precisely those terms for which
the cycles $a$ and $b$ have mod-p linking number zero with $U$ will
contribute. For a given $(\epsilon_1,\epsilon_2)$, the linking
numbers will be (giving $U$ the clockwise orientation):
\begin{eqnarray*}
lk(a,U) & = & \epsilon_1, \\
lk(b,U) & = & \epsilon_2 + 1.
\end{eqnarray*}
If we collect this information into the monomial 
$F_{\epsilon_1,\epsilon_2}(t_a,t_b) 
= t_{a}^{lk(a,U)} t_{b}^{lk(b,U)}$,
then for a given $p$, for some determinable
sign $\alpha$:
\[
\lambda_{CWL} = \alpha 2 p \sum_{(\epsilon_1,\epsilon_2)=(0,0)}^{(1,1)}
\left(
\frac{1}{p}\sum_{r=0}^{p-1} \frac{1}{p}\sum_{s=0}^{p-1} F_{\epsilon_1,\epsilon_2}( e^{r \frac{2\pi i}{p}},
e^{s \frac{2\pi i}{p}} )
\right)
.
\]

[ To see this, note the following identity:
\[
\sum_{s=0}^{p-1} ( e^{\frac{2\pi i}{p}s}  )^t = \left\{
\begin{array}{ll}
p & \mbox{if}\ p|t, \\
0 & \mbox{otherwise.]}
\end{array}
\right.
\]

We can seperate out the dependence on $p$ as follows. Defining the polynomial:
\[
F'(t_a,t_b) = \sum_{(\epsilon_1,\epsilon_2) = (0,0)}^{(1,1)} F_{\epsilon_1,\epsilon_2}(t_a,t_b),
\]
then the result is:
\[
\lambda_{CWL}( \Sigma^p_{U^{\kappa}} ) = 
\alpha \frac{2}{p} \sum_{r=0}^{p-1} \sum_{s=0}^{p-1}
F'(e^{r \frac{2\pi i}{p}}, e^{s \frac{2\pi i}{p}} ).
\]

\end{example}

\begin{remark}
Observe that such a procedure can be implemented for any decoration of 
any knot by a complete graph Y-link whose underlying diagram is a theta
graph with some legs attached. Thus we have a flexible new 
construction of knots for which the Casson-Walker-Lescop invariants
of their $p$-fold branched cyclic covers may be calculated, for all choices
of $p$. It is interesting that this sequence is determined by the values of
an associated polynomial at the roots of unity, evoking the familiar
formula for the order of the first homology in terms of the Alexander-Conway
polynomial \cite{F,HK} of the branching knot.

Observe also that similar constructions and 
arguments can be applied to any projection of the
LMO invariant (see Example \ref{examp}).

\end{remark}

Section \ref{bcc} will prove part 1 of the main theorem.
Section \ref{exprsect} will prove the Lemmas \ref{calclemma} and \ref{calclemmab} which facilitate the extension.

\subsection{Branched cyclic covers}\label{bcc}

Excellent expositions of this construction 
can be found in the literature, 
for example in \cite{Rol} and in \cite{Lic}.

Let $N(K)$ be a regular neighbourhood of a knot $K$ in an integral
homology three-sphere $M$, and let
$X_K$ be the closure of its exterior. Let ${X^p_K}$ be 
the $p$-fold cyclic cover of $X_K$.


In this work, it will be convenient to construct this space directly
with a Seifert surface $F$ for $K$. 
This direct construction uses
the space $W_F$, which is obtained by  ``splitting $X_K$
open along $F$''. In practice, $W_F$ is obtained by removing the intersection 
of an open collar neighbourhood of the surface $F$ with $X_K$: that is,
if $N: F\times [-1,1] \hookrightarrow M$
is a bicollar for $F$ then $W$ is obtained as the removal of $X_K\cap N(F\times
(-\epsilon,\epsilon))$ from $X_K$. The subspaces $X_K\cap N(F\times \{\epsilon\})$ and $X_K\cap N(F\times \{-\epsilon\})$ are  
homeomorphic images of $F$, call them $F^+$ and $F^-$, and there will
exist a homeomorphism between them $\phi: F^- \rightarrow F^+$, a gluing 
along which recovers the exterior of the knot, $X_K$.

The unique
$p$-fold cyclic cover of the exterior of 
the knot, ${X^p_K}$, is constructed from 
$p$ copies of $W_F$: denote them $W_F^i$, where the index $i$ is an element
of ${\Zset_p}$. 
Glue $F^-$ in $W_F^i$ to $F^+$ in $W_F^{i+1}$ using $\phi$ composed 
with the identification $W_F^i \simeq W_F^{i+1}$.

The $p$-fold cyclic cover of $M$ {\it branched} over $K$ is obtained
as a completion of the space ${X^p_K}$. Specifically, a solid torus
is glued into the torus boundary of ${X^p_K}$
so that its meridian projects to a meridian of $K$ traversed $p$ times.

In this work, we will encounter a situation where a knot $K^L$ is recovered
from another knot $K$ by surgery on some framed link $L$ decorating $K$, 
and where there is given some Seifert surface $F$ for $K$ in the complement of
$L$. 
In such a situation, a Seifert surface $F^L$ for $K^L$ may be canonically 
chosen: namely, just choose the position of $F$ in the surgered manifold.
The 
pair $(W_{F^L},\phi^L)$ can then be obtained from the pair $(W_F,\phi)$ by
surgery on the link $L$ in $W_F$. Thus, by construction, $\Sigma_{K^L}^p$
is recovered by surgery on the
framed link in $\Sigma_K^p$ comprising of one copy of the original framed
link in each of the $p$ copies of $W_F$ included in $\Sigma_K^p$.

\subsection{Complete graph Y-links}
\label{manoeuvres}

Before we describe this proof, we introduce a term.

\begin{definition}
A {\bf mixed Y-link} in a three-manifold
is a an embedding of a set of Y-components {\it and possibly}
some claspers, into that three-manifold.
\end{definition}

Let us turn to the scenario described by the main theorem of this section.
We have two knots: $K$ and a knot $K^L$ obtained from $K$ by surgery on a 
complete graph Y-link $L$ of surplus $s$. 
The strategy is to replace $L$ with a mixed Y-link also presenting
the knot $K^L$, but which is in the complement of a Seifert surface
$\check{F}$ for $K$. 

To begin, choose some Seifert surface $F$. We will
now isotope the graph Y-link $L$ so that it is in a special position with
respect to $F$. We require two conditions of this position. The first 
requirement is that the band of a clasper associated to a leg does not 
intersect $F$ (although one of the leaves certainly will). 
Given such a position, which can always be obtained,
then around every leaf we can find a ball so that the arrangement is
obtained by some homeomorphism of that ball with the following ball
in $S^3$:

\vspace{0.6cm}

\[
\leg
\]\vspace{2.3cm}

The second requirement is that the only other intersections of $L$ with $F$
are transversal intersections of bands with $F$. 

A graph Y-link in such a position will now be modified in two steps. 
The first step is make an
adjustment of $L$
around each leg, where it appears as in the ball above. The adjustment
is made with the move described in Lemma \ref{legbreak}:\vspace{1.8cm}

\[
\legged
\]\vspace{2.9cm}

There is one such adjustment to be made for every (if there are any, that is)
leg of the graph Y-link, leaving the decoration with 
one transverse intersection of a band with 
$F$ for every leg, and with possible further transverse 
intersections of internal bands with $F$.

The second step is to adjust the resulting
mixed graph Y-link {\it and Seifert surface}
at each such intersection (so
one for each leg with possible further adjustments).
Using the move described in Lemma \ref{insertclasp} the band can be broken into
a clasp
at that intersection,
then $F$ can be tubed
around the introduced annulus (this is shown in the diagram 
following the proof).

After these adjustments, 
we have a mixed graph Y-link in the complement of a Seifert surface
for $K$. 
In this mixed graph Y-link, one can identify $s$ Y-components (one for
each surplus trivalent vertex) and a possible number of other claspers
that have emerged from the modifications. 
By an abuse of notation call this mixed Y-link $L$. 
This has been constructed to be in the complement of some
Seifert surface for $K$, call it $\check{F}$.

Use this Seifert surface to construct $\Sigma^p_K$.
Denote by $\widetilde{L}$ the
Y-link comprised of a copy of $L$ for each of the
$p$ copies of $W_{\check{F}}$ included into $\Sigma^p_K$. 
By construction, surgery on this mixed Y-link 
yields $\Sigma^p_{K^L}$.

\

\vspace{1cm}

\

\[
\pierce\hspace{2cm} \rightarrow \hspace{3.5cm}
\Picture{
\put(0,2){$\modhook$}}
\]\vspace{2cm}

\subsection{Proof of Theorem \ref{maintha}}

\

Our task is to show that $\Sigma^p_K$ and $\Sigma^p_{K^L}$ are $s-1$-equivalent.
It is sufficient, by Lemma \ref{nequivlem}, to exhibit an  $s$-scheme 
relating $\Sigma^p_{K^L}$ to $\Sigma^p_K$. 
This $s$-scheme will be
based upon the mixed Y-link $\widetilde{L}$ in $\Sigma^p_{K}$.

The appropriate selection of Y-sublinks is immediate. 
Order the Y-components of $L$, denoting them
Y$_1$ through to Y$_s$. 
Let $\widetilde{L}_i$ denote the Y-sublink of $\widetilde{L}$
comprised of the set of $p$ copies of Y$_i$.
These $s$ disjoint Y-sublinks form the required 
relating $s$-scheme: $\{\widetilde{L}; \widetilde{L}_1,\ldots,\widetilde{L}_s\}$.
Recall that $\widetilde{L}_{i_1,\ldots,i_s}$  denotes the link obtained
by forgetting those sublinks $\widetilde{L}_{j}$ for which $i_j=1$.
By construction, the result of 
surgery on $\widetilde{L}_{i_1,\ldots,i_s}$
is the $p$-fold branched cyclic cover
branched over the knot obtained from $K$
by doing surgery on $L_{i_1,\ldots,i_s}$. If this multiplet is not
$\{0,\ldots,0\}$ then this knot will be $K$.
\rtb

\

\subsection{A convenient expression}

Lemma \ref{niceexprlem} can be employed to give an expression 
for $\Sigma^p_{K^L} - \Sigma^p_K$ that will be employed frequently
in the sequel. Let $\widetilde{L}^{(a_1,\ldots,a_s)}$ denote the
sublink of $\widetilde{L}$ that forgets, for each $i$ between $1$
and $s$, every copy of Y$_i$ except the
copy in the included subspace $W_{\check{F}}^{a_i}$.
\begin{lemma}\label{nelt}
\[
\Sigma^p_{K^L} - \Sigma^p_{K} = 
\sum_{(a_1,\,\ldots,a_s)=(1,\,\ldots,1)}^{(p,\,\ldots,p)}
[ \Sigma^p_K , \widetilde{L}^{(a_1,\,\ldots,a_s)} ].
\]
\end{lemma}

\section{Proofs of the lemmas}\label{exprsect}

Before we prove these lemmas, we must first recollect
how to compare two graph Y-links in a given three-manifold.

\subsection{Comparing graph Y-links in three-manifolds} \label{orientquest}
\label{signsect}

To an n-component graph Y-link $L$ in a three-manifold $M$ 
is canonically associated 
a trivalent graph $D_L$: take a trivalent vertex
for every Y-component, and join tines with an edge
when the associated leaves are
linked in $M$. Note that, at this point, we are not associating any 
extra structure to this diagram (like cyclic orderings of edges
at trivalent vertices, or orientations of edges). 
Let $L'$ be another $n$-component
graph Y-link and let there exist a graph isomorphism
between their associated diagrams $D_L$ and $D_{L'}$.
\begin{lemma}\label{complem}
\[
[M,L] = \pm [M,L'] \in \CG_n^Y\CM.
\]
\end{lemma}
\underline{Proof.} Take some surgery link presenting $M$. 
The clasper graphs presenting 
the Y-links
$L$ and $L'$ can be be isotoped in $M$
to appear in the complement of the surgery
presentation in $S^3$. We now ``move'' the position of
$L'$ to that of $L$. That is, first isotope the vertices (discs of the
clasper graph) of $L'$ 
(in the complement of the surgery link) to occupy the same position
as the corresponding vertices of $L$, choosing the identification
of two discs which aligns the appropriate edges.
Now $L'$ can be adjusted to be in the position of
$L$ by crossing changes between its
internal bands, and by crossing changes between its internal bands and surgery
components (using Lemma \ref{crosstube}). A number of half-twists will 
possibly be introduced as the
last step (using Lemma \ref{minusintro}).
\rtb

In this work, we will employ the following procedure to determine
the sign of the above comparison.
Begin by choosing an orientation for each
Y-component of $L$ (that is, orient the associated thrice-punctured disc).
This induces an orientation on the associated boundary components, which
will appear as follows (or its mirror image):
\vspace{2cm}

\hspace{5.8cm}\actycompa
\vspace{2.6cm}

A graph Y-link with such a choice at every Y-component
will be referred to as an {\it oriented graph Y-link}, and the choice
will be called an {\it orientation} for it.
This choice is equivalent to a choice of cyclic ordering of the 
leaves at every Y-component. To fix this correspondence, for 
example, say that the edges of the above Y-component have a counter-clockwise
cyclic ordering. Record this information on the associated
graph (that is, cyclically order the edges at every trivalent vertex
of it). Call this enhancement an {\it orientation} for the
graph.

To each edge of an oriented graph, associated to an oriented
graph Y-link, a sign will be associated: the {\it twist} of that edge. 
The twist is defined as the linking number of the (oriented) inner
circles of the two leaves whose linking that edge represents.
Let $T(e)$ denote the twist of the edge $e$.

This information can now be used to compare two graph Y-links, $L$ and $L'$,
equipped with an isomorphism $\phi: D_L \rightarrow D_{L'}$. 
Namely, choose an orientation
for $L$, inducing an orientation on $D_L$. The graph isomorphism induces
an orientation on $D_{L'}$, and hence on $L'$.
\begin{lemma}\label{enhanced}
\[
[M,L] = \left(\,\prod_{\mbox{$e$ an edge}} {T(e)}{T'(\phi(e))}\right)
[M,L'] \in \CG^Y_n\CM
\]
\end{lemma}
\underline{Proof.}
We enhance the proof of Lemma \ref{complem}. Notice that the crossing changes
that are made as $L'$ is ``moved'' to $L$ do not affect the 
twist of any edge. Introducing a half twist, however, affects the twist of the 
associated edge by multiplying it by minus one.
\rtb

It will also prove necessary, in the sequel, to have an expression
for the twist of an edge whose associated  
pair of leaves are
linked indirectly by a chain of claspers. 
That is, where a leaf $L_i$ of a Y-component $Y_i$ is linked in a ball with
a leaf of a clasper $C_1$, whose other leaf links in a ball with
a leaf of another clasper $C_2$, and so on, until the clasper $C_\mu$ whose
other leaf links a leaf $L_f$ of the other Y-component $Y_f$.

To calculate the twist in such a situation (that is, the twist that 
results from surgeries on all those claspers),
begin by choosing an orientation for each clasper in the chain
(that is, an orientation for the twice-punctured disc). It will
appear as below (or its mirror image):
\[
\setlength{\unitlength}{30pt}
\hspace{1.3cm}
\Picture{
\put(-1.4,-0.3){\vector(0,-1){0.01}}
\put(-1.2,-0.1){\vector(0,1){0.01}}
\put(1.2,-0.3){\vector(0,-1){0.01}}
\put(1.4,-0.1){\vector(0,1){0.01}}
\qbezier(1.1,0.1)(1.4,0.1)(1.4,-0.2)
\qbezier(1.1,-0.5)(1.4,-0.5)(1.4,-0.2)
\qbezier(1.1,-0.1)(1.2,-0.1)(1.2,-0.2)
\qbezier(1.1,-0.3)(1.2,-0.3)(1.2,-0.2)
\qbezier(1.1,-0.5)(0.9,-0.5)(0.9,-0.5)
\qbezier(1.1,-0.3)(0.9,-0.3)(0.9,-0.3)
\qbezier(0.9,-0.5)(0.6,-0.5)(0.6,-0.2)
\qbezier(0.9,-0.3)(0.8,-0.3)(0.8,-0.2)
\qbezier(0.9,-0.1)(0.8,-0.1)(0.8,-0.2)
\qbezier(0.9,0.1)(0.6,0.1)(0.6,-0.2)
\qbezier(-1.1,-0.5)(-1.4,-0.5)(-1.4,-0.2)
\qbezier(-1.1,0.1)(-1.4,0.1)(-1.4,-0.2)
\qbezier(-1.1,-0.3)(-1.2,-0.3)(-1.2,-0.2)
\qbezier(-1.1,-0.1)(-1.2,-0.1)(-1.2,-0.2)
\qbezier(-1.1,0.1)(-0.9,0.1)(-0.9,0.1)
\qbezier(-1.1,-0.1)(-0.9,-0.1)(-0.9,-0.1)
\qbezier(-0.9,0.1)(-0.6,0.1)(-0.6,-0.2)
\qbezier(-0.9,-0.1)(-0.8,-0.1)(-0.8,-0.2)
\qbezier(-0.9,-0.3)(-0.8,-0.3)(-0.8,-0.2)
\qbezier(-0.9,-0.5)(-0.6,-0.5)(-0.6,-0.2)
\qbezier(0.6,-0.1)(-0,-0.1)(-0.6,-0.1)
\qbezier(0.6,-0.3)(0,-0.3)(-0.6,-0.3)
\qbezier(-1.1,0.1)(-1.4,0.1)(-1.4,-0.2)
\qbezier(-1.1,-0.5)(-1.4,-0.5)(-1.4,-0.2)
\qbezier(-1.1,-0.1)(-1.2,-0.1)(-1.2,-0.2)
\qbezier(-1.1,-0.3)(-1.2,-0.3)(-1.2,-0.2)
\qbezier(-1.1,-0.5)(-0.9,-0.5)(-0.9,-0.5)
\qbezier(-1.1,-0.3)(-0.9,-0.3)(-0.9,-0.3)
\qbezier(-0.9,-0.5)(-0.6,-0.5)(-0.6,-0.2)
\qbezier(-0.9,-0.3)(-0.8,-0.3)(-0.8,-0.2)
\qbezier(-0.9,-0.1)(-0.8,-0.1)(-0.8,-0.2)
\qbezier(-0.9,0.1)(-0.6,0.1)(-0.6,-0.2)
\qbezier(--1.1,-0.5)(--1.4,-0.5)(--1.4,-0.2)
\qbezier(--1.1,0.1)(--1.4,0.1)(--1.4,-0.2)
\qbezier(--1.1,-0.3)(--1.2,-0.3)(--1.2,-0.2)
\qbezier(--1.1,-0.1)(--1.2,-0.1)(--1.2,-0.2)
\qbezier(--1.1,0.1)(--0.9,0.1)(--0.9,0.1)
\qbezier(--1.1,-0.1)(--0.9,-0.1)(--0.9,-0.1)
\qbezier(--0.9,0.1)(--0.6,0.1)(--0.6,-0.2)
\qbezier(--0.9,-0.1)(--0.8,-0.1)(--0.8,-0.2)
\qbezier(--0.9,-0.3)(--0.8,-0.3)(--0.8,-0.2)
\qbezier(--0.9,-0.5)(--0.6,-0.5)(--0.6,-0.2)
\qbezier(-0.6,-0.1)(--0,-0.1)(--0.6,-0.1)
\put(0,-0.3){\vector(1,0){0.01}}
\put(0,-0.1){\vector(-1,0){0.01}}
\qbezier(-0.6,-0.3)(-0,-0.3)(--0.6,-0.3)
}
\]\vspace{1cm}

Let $l_0$ be the linking number of the inner circles of the annuli
corresponding to the appropriate pair of leaves, one from $L_i$ and
one from $C_1$; let
$l_1$ be the linking number of the inner circles of the 
appropriate pair of leaves, one from $C_1$ and one from
$C_2$; and so on. The following fact is straightforward:

\begin{lemma}\label{twistcalc}
The following product 
is well-defined (independent of the choices made
in orienting the claspers) and 
\[
T(e) = l_0\prod_{i=0}^{\mu}(-l_i).
\]
\end{lemma}
\begin{remark}
Actually, this is more information than we will use. All we need to know
is that the resulting twist is a function of the linking numbers along the
chain. One can prove this fact by considering the effect on the left-most 
linking number of a surgery on the left-most clasper; and proceed inductively.
\end{remark}

\subsection{Proof of Lemma \ref{welldefnlem}}

\begin{lemma}\label{welldefnlem}
$(\Sigma^p_K)^{\lambda_L'}
- \Sigma^p_K$ is well-defined in $\CG_s^Y\CM$.
\end{lemma}
\underline{Proof.}
If the Y-components of $L$ are oriented (according to the Section \ref{orientquest}), 
then any graph Y-links in $\Sigma^p_K$ resulting from the definition
of $\lambda_L'$
will have isomorphic
underlying vertex-oriented graphs with identical twists along the edges, so
that the well-definedness follows from Lemma \ref{enhanced}.
\rtb

\subsection{Proof of Lemma \ref{calclemma}}

\

The graph Y-link $L$ is modified in the ball, following
the manoeuvres of Section \ref{manoeuvres}, to appear in that
ball in the following way: 

\

\vspace{3cm}

\hspace{6cm}\leggedp\vspace{5.5cm}

Surgery on this mixed Y-link recovers the knot $K^L$. The knot $K^{L'}$
is recovered by surgery on the sublink that is obtained by forgetting
the claspers B and C, and the knot $K^{L''}$ is obtained by surgery
on the subgraph that is obtained by forgetting the clasper A and smoothing
the half-twist. By an abuse of notation, we will 
refer to these
mixed Y-links as $L$, $L'$ and $L''$.

The branched cyclic covering $\Sigma^p_{K^L}$ (resp. $\Sigma^p_{K^{L'}}$,
resp.
$\Sigma^p_{K^{L''}}$) is obtained from $\Sigma^p_{K}$ by performing
surgery on the mixed Y-link that is comprised of a copy of
$L$ (resp. ${L'}$, resp. ${L''}$) for each of the
$p$ subspaces $W_{\check{F}}$ included in $\Sigma^p_K$. Call these mixed
Y-links $\widetilde{L}$,
$\widetilde{L'}$ and $\widetilde{L''}$.
The mixed Y-link $\widetilde{L}$ will have $p$ copies of each of the components
illustrated above. This set of components will be denoted by affixing
a $\Zset_p$-valued superscript 
to the labels Y,A,B and C. 
We will use the same symbols to denote the corresponding components in the
sublinks $\widetilde{L'}$ and $\widetilde{L''}$.

To fix this notation fully, it remains
to specify which
boundary component after the removal of $\check{F}$ corresponds to $\check{F}^+$. This choice can be specified 
if we declare that we are making the choice
that results in $C^i$ linking $B^{i+1}$ in a ball, as follows:

\

\vspace{1.7cm}

\hspace{5.7cm}\seifcross
\vspace{2cm}

Order the Y-components of $L$ so that the one pictured above is Y$_1$.
Let $\widetilde{L}^{(a_1,\ldots,a_s)}$ denote the sublink of 
$\widetilde{L}$ that is obtained by forgetting every Y-component, except
the copy of Y$_i$  in $W_{\check{F}}^{a_i}$, for each $i$. 
Similarly use the symbol $\widetilde{L'}^{(a_1,\ldots,a_s)}$
(resp. $\widetilde{L''}^{(a_1,\ldots,a_s)}$) for the sublink obtained 
from $\widetilde{L}^{(a_1,\ldots,a_s)}$ by forgetting every
B$^i$ and C$^j$ (resp. forgetting every A$^k$ and smoothing the
half-twist in every B$^{l}$).

Now, Lemma \ref{nelt} indicates that
\[
\Sigma^p_{K^L} - \Sigma^p_{K} = 
\sum_{(a_1,\ldots,a_s)=(1,\ldots,1)}^{(p,\ldots,p)}
[ \Sigma^p_K , \widetilde{L}^{(a_1,\ldots,a_s)} ],
\]
with identical expressions holding when $L$ is replaced by $L'$ or
$L''$. The lemma in question 
then follows from the following equation, which we will subsequently 
demonstrate:
\begin{equation}\label{cruxeqn}
[ \Sigma^p_K , \widetilde{L}^{(a_1,\ldots,a_s)} ] =
[ \Sigma^p_K , \widetilde{L'}^{(a_1,\ldots,a_s)} ] -
[ \Sigma^p_K , \widetilde{L''}^{(a_1,\ldots,a_s)} ].
\end{equation}

Consider the mixed Y-link
$\widetilde{L}^{(a_1,\ldots,a_s)}$ corresponding to some
$s$-tuplet $(a_1,\ldots,a_s)$. Let the mixed Y-links 
$R^{(a_1,\ldots,a_s)}$ and $S^{(a_1,\ldots,a_s)}$
be obtained from it by modifying it in the included subspace $W_{\check{F}}^{a_1}$
as follows:

\

\vspace{2.25cm}

\hspace{6cm}\tripa\hspace{0.5cm},

\

\vspace{2cm}

\hspace{3.3cm}\tripb\hspace{1.25cm}and\hspace{4cm}\tripc\hspace{0.5cm}.

\vspace{1.2cm}

Lemma \ref{add}, indicates that
\[
[ \Sigma^p_K , \widetilde{L}^{(a_1,\ldots,a_s)} ] =
[ \Sigma^p_K , {R}^{(a_1,\ldots,a_s)} ] +
[ \Sigma^p_K , {S}^{(a_1,\ldots,a_s)} ].
\]
Equation \ref{cruxeqn} follows from the identifications
\begin{eqnarray} 
\ [ \Sigma^p_K , {R}^{(a_1,\ldots,a_s)} ] & = &-\ [ \Sigma^p_K , \widetilde{L''}^{(a_1,\ldots,a_s)} ], \\
\ [ \Sigma^p_K , {S}^{(a_1,\ldots,a_s)} ] & = &\ \ \ 
[ \Sigma^p_K , \widetilde{L'}^{(a_1,\ldots,a_s)} ].
\end{eqnarray}

We leave the reader to check these relations. The identifications proceed
by removing claspers when they have a leaf that bounds a disk in the complement of the rest of the graph. Note that the half-twist in B$^{a_1}$ 
is smoothed at the
expense of the minus one
(using Lemma \ref{minusintro}).
\rtb

\subsection{Proof of Lemma \ref{calclemmab}}\label{calcbsect}


Choose a Seifert surface $F$ for $K$, and an orientation for $K$, and hence
an orientation for the Seifert surface. Represent this orientation by
a normal vector field. Denote the embedded trivalent graph 
that is associated to a position of $L$ that only meets this Seifert
surface in transversal intersections of internal edges, by $E_L$.

Let $L$ be used again (abusing the notation) 
to denote the mixed Y-link that results from the
modifications of this embedding with respect to this Seifert surface
that are described in  Section \ref{manoeuvres}.
In this case, there are no legs to worry about, so the only modifications 
occur when edges of the associated clasper graph intersect the Seifert surface.
\vspace{1cm}

\[
\pierce\hspace{2cm} \rightarrow \hspace{3.5cm}
\Picture{
\put(0,2){$\modhook$}}
\]\vspace{1.6cm}

The mixed Y-link $L$ comprises of $s$ Y-components
(order them Y$_1,\ldots,$Y$_s$)
and a number of other claspers, in the complement of some Seifert surface
$\check{F}$. 




Recall that Lemma 
\ref{nelt}
indicates that
\[
\Sigma^p_{K^L} - \Sigma^p_{K} = 
\sum_{(a_1,\ldots,a_s)=(1,\ldots,1)}^{(p,\ldots,p)}
[ \Sigma^p_K , \widetilde{L}^{(a_1,\ldots,a_s)} ],
\]
using notation explained there.

Let us consider one of these terms, corresponding to some 
$p$-tuplet $(a_1,\ldots,a_s)$. We will show 
that corresponding to each edge of $E_L$,
to whose endpoints the two Y-components Y$_i$ and Y$_j$ are associated,
there will be an equation in $\Zset_p$ relating 
$a_i$ and $a_j$ which must be satisfied in order
that the corresponding term $[ \Sigma^p_K , \widetilde{L}^{(a_1,\ldots,a_s)} ]$
be non-zero.

Orient the edges of $E_L$. Define functions $i$ and $f$ from the edges
of $E_L$ to the labels of the Y-components of $L$ so that $Y_{i(e)}$ 
corresponds to the origin of the edge $e$, and the Y-component $Y_{f(e)}$
corresponds to its end. 
Let $<e,F>$ denote the signed sum of intersections of the edge $e$
(oriented from $i(e)$ to $f(e)$) 
with $F$ (where a plus is an intersection in the direction of the 
orienting normal vector field). 
\begin{lemma}

\

If for the $s$-tuplet $(a_1,\ldots,a_s)$, 
\begin{equation}
a_{f(e)} \neq a_{i(e)} + <e,F>\ \in \Zset_p,
\end{equation}
for some edge $e$ of $E_L$, then
\[
[ \Sigma^p_K , \widetilde{L}^{(a_1,\ldots,a_s)} ] = 0.
\]
\end{lemma}
\underline{Proof}

\
 
Take such a pair of Y-components in $L$, Y$_{i(e)}$ and Y$_{f(e)}$. Every time
the corresponding edge of $E_L$ intersected the Seifert surface $F$ it
was broken into a clasp. Thus, in general, 
Y$_{i(e)}$ and Y$_{f(e)}$ will be linked in $L$ via some chain of claspers. 
Denote them 
$C_1$ up to $C_\mu$ where a leaf of $C_1$ links a leaf of Y$_{i(e)}$, with
its other leaf linking a leaf of the clasper $C_2$ etc.
Denote the copies of each of these claspers that occur
in the graph $\widetilde{L}^{(a_1,\ldots,a_s)}$ by affixing
a $\Zset_p$-valued 
superscript.

To show this equation, 
we proceed to remove copies of the $C_j$ when they have a leaf that bounds
a disc 
(in the complement of the rest of the mixed Y-link is to be understood
when we use this phrase). 
If, finally, we are left with a 
mixed Y-link where one
of the Y-components has a leaf that bounds a disk, then we know that
the contribution $[ \Sigma^p_K , \widetilde{L}^{(a_1,\ldots,a_s)} ]$
will be zero.

If the edge $e$
has no 
(resp. 1 positive, resp. 1 negative) intersection 
with the Seifert surface $F$, then the
term
$[\Sigma^p_K,\widetilde{L}^{(a_1,\ldots,a_s)}]$ will be zero unless
$a_{f(e)} = a_{i(e)}$ (resp. $a_{f(e)} = a_{i(e)}+1$,
resp. $a_{f(e)} = a_{i(e)}-1$), for
otherwise the Y-components will have leafs bounding discs.

Consider, then, the situation
where the edge $e$ does intersect $F$ at at least two points, so that
some extra claspers arise. 
If the first intersection 
is a positive (resp. negative) 
intersection, then we can remove all copies of $C_1$ 
except the copy $C_1^{a_{i(e)}+1}$ (resp. $C_1^{a_{i(e)}-1}$). 
Proceeding, we remove all copies of $C_2$ except the one that intersects
the surviving copy of $C_1$. And so on. One can check that as a result
of this process, the surviving copy of $C_\mu$ will intersect precisely
the Y-component $Y_j^{a_{i(e)}+<e,F>}$. 
\rtb

\begin{lemma}

\

If the $s$-tuplet $(a_1,\ldots,a_s)$ 
satisfies the equations
\begin{equation}\label{lifteqn}
a_{f(e)} = a_{i(e)} + <e,F>\ \in \Zset_p,
\end{equation}
for every edge $e$ of $E_L$, then
\[
[ \Sigma^p_K , \widetilde{L}^{(a_1,\ldots,a_s)} ] = (\Sigma^p_K)^{\lambda_L'}
- \Sigma^p_K\ \ \in \CG^Y_s\CM.
\]
\end{lemma}
\underline{Proof.}

\

Take the mixed Y-link $L$ that is obtained following Section 
\ref{manoeuvres}, where every intersection with a Seifert surface
is broken into a clasp. Orient each Y-component and clasper of $L$.
This introduces an orientation on each Y-component and clasper of 
$\widetilde{L}^{(a_1,\ldots,a_s)}$. Forget every clasper with a leaf
that bounds a disc (as in the previous proof).

Note that the linking numbers
between pairs of bounding circles in $L$ and the corresponding
pairs in $\widetilde{L}^{(a_1,\ldots,a_s)}$ are the same. 
Thus, the theorem follows from 
Lemma \ref{enhanced} and Lemma \ref{twistcalc}. 
\rtb

It is clear that the Equations \ref{lifteqn} are precisely the requirement
that every loop on $E_L$ have zero mod-$p$ linking number with $K$.
Note that when at least one solution exists, there are precisely
$p$ solutions, corresponding to the possible values of $a_1$.
\rtb

\

\section{Part 2: The relation to Vassiliev theory}
\label{realproof}

In Part 2 we consider some very general questions about the relationship
with Vassiliev theory of 
knot invariants obtained by composing finite-type three-manifold
invariants with the branched cyclic covering construction.

The Vassiliev theory of finite type
invariants of oriented knots in the three-sphere centers on a filtration
of $\CK(S^3)$, the free abelian group generated by ambient isotopy 
classes of oriented knots:
\[
\CK(S^3) \supset  \CF_1\CK \supset \CF_2\CK \supset \ldots
\]

Our first speculation might be that any finite-type three-manifold invariant
$v$, composed with $\Sigma^p$, is a finite-type invariant of knots. To
begin, then, we will see that this hypothesis can be ruled out using
simple arguments (for any choice of $p$). 
The method of proof also provides some other interesting 
information: that the function $|H_1(\Sigma^p_K)|$ can be unbounded on a
sequence of knots of strictly increasing $n$-triviality.

\begin{theorem}

Fix a positive integer $p$. For every integer $n>0$, there exists
a non-trivial element $E_n \in \CF_n\CK$ such that
\begin{equation}
\Sigma^p(E_n) \notin \CF^Y_1\CM.
\end{equation}
\end{theorem}

\underline{Proof.}

\

Let $\Omega_{2}$ be the knot obtained from surgery on the 
following decoration of the unkot by a 0-framed surgery link (that is,
the link obtained by replacing twice-punctured discs by surgery
pairs, as in Appendix \ref{GHT}). Similarly
let $\Omega_{n}$ denote the obvious extension to a knot obtained from
surgery on a wheel with $n$ spokes.
\begin{equation}
\epsfxsize = 4cm
{\fig{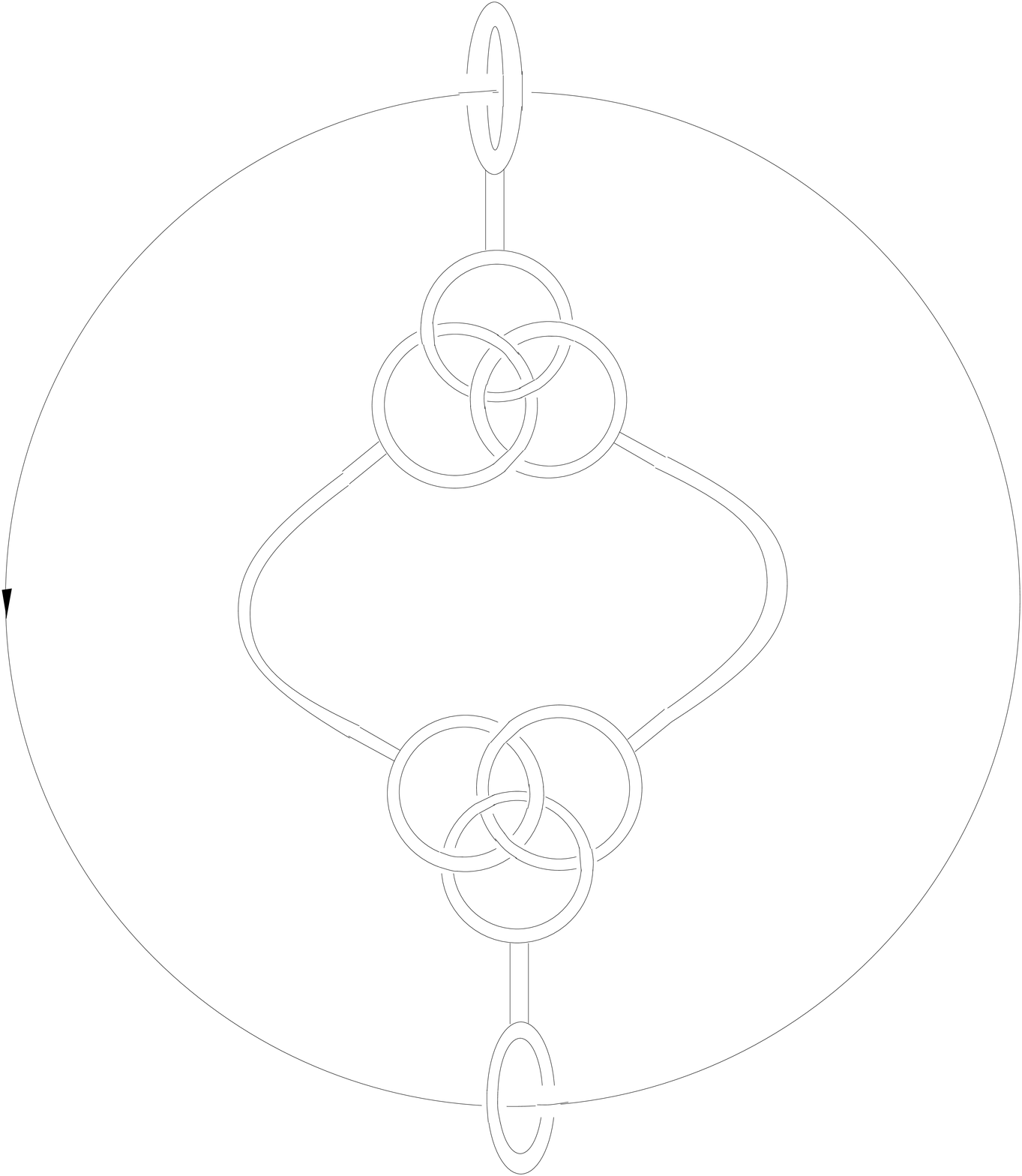}}\ \ .
\end{equation}
$\Omega_n$ is $(n-1)$-trivial, being obtained from surgery
on a graph Y-link of degree $n$ (alternatively we can identify it as a 
Kanenobu-Ng wheel \cite{Kan,Ng}).

Our theorem will follow if we can identify an $n'\geq n$ so that 
$\Sigma^p(\Omega_{n'})$ is not an integral homology three sphere. Then
the proof would be completed by setting $E_{n}=\Omega_{n'}-U$, where
$U$ is the unknot (see
Matveev's theorem, Theorem \ref{Mattheo}).

To this end we can exploit a formula for the order of the first homology
of the $p$-fold branched cyclic covering of a knot $K$ in terms of its
Alexander polynomial $A_K(t)$ \cite{G,F,HK}. (The symbol $|H_1( \Sigma^p(K) ) |$ below is extended to be zero if $\Sigma^p(K)$ has non-vanishing first Betti
number).
\[
| H_1( \Sigma^p(K) ) | = \prod_{q=0}^{p-1} | A_K( e^{2 \pi i \frac{q}{p}} ) |.
\]
In \cite{K} we calculated (alternatively, (do the work to) identify 
this as one of the knots 
considered in \cite{Kan})
\[
A_{\Omega_n}(t) = (1-(1-t)^{n})(1-(1-t^{-1})^{n}),
\]
so that
\[
| H_1 ( \Sigma^p(\Omega_n) )| = \prod_{q=0}^{p-1} |( 1 - (1- e^{2\pi i \frac{q}{p}})^n)|^2.\]

Define a function $f(p,n)$ to be $|H_1 ( \Sigma^p(\Omega_n) )|$. 
We are interested in the behaviour of $f(p,n)$ as $n$ goes to infinity whilst
$p$ is fixed.

Consider a factor $(1 - \alpha(p,q)^{n})$ in the above, where 
$\alpha(p,q) = 1- e^{2\pi i \frac{q}{p}}$. If $|\alpha(p,q)| < 1$, then
$\lim_{n\rightarrow \infty}|1 - \alpha(p,q)^{n}| = 1$. Alternatively,
if $|\alpha(p,q)| > 1$, then
$\lim_{n\rightarrow \infty}|1 - \alpha(p,q)^{n}| = \infty$.
Note, however, that there are precisely two possible choices of $\frac{q}{p}$
such that $|\alpha(p,q)| = 1$; namely $\frac{q}{p}=\pm \frac{1}{6}$ 
(which will occur in $f(p,q)$ whenever $6|p$). In these cases
$\alpha(6,\pm 1) = e^{\mp 2\pi i \frac{1}{6}}$, and  
$\lim_{n\rightarrow \infty}(1 - \alpha(6,\pm 1)^{n})$ does not exist.
We will get around this by choosing a subsequence; namely
$\lim_{n\rightarrow \infty}(1 - \alpha(6,\pm 1)^{6n+3}) = 2.$

Finally, observe that there is always at least one factor $|1-\alpha(p,q)|$ 
such that $|\alpha(p,q)|>1$. Then
\[
\lim_{n\rightarrow \infty} f(p,6n+3) = \infty.
\]
\rtb
See \cite{Gor} for similar techniques.

We next consider a measure of the independence of these theories.
The following theorem is modelled on Stanford--Trapp's
definition of an invariant which is 
``independent of finite-type invariants'',
although note that the stated property (appears to be) 
substantially weaker.

It says that we may realise (some vector proportional to) any combination
of primitive diagrams as the first non-vanishing term of the LMO
invariant on the
$p$-fold branched cyclic covering branched over a knot which is $n$-trivial, 
{\it where the $n$-triviality is a free parameter}.

According to the notation of \cite{LeGr}, $\hat{Z}^{LMO}$ 
denotes the LMO invariant normalised by powers of the order of the first 
homology (the normalisation that satisfies the simple connect-sum formula). 

\begin{theorem}\label{notft}
Take positive integers $n,p$ and $s$, and a knot $K$ in $S^3$
such that the $p$-fold branched 
cyclic covering of $S^3$ branched over $K$
is a rational homology three-sphere. 
Then, there exists a non-zero integer $\alpha(s,n,p)$
such that for every 
$\Zset$-linear combination of primitive degree $s$
uni-trivalent diagrams on an empty skeleton, $D$, there exists a knot
$K(s,n,p,D)$ satisfying:
\begin{itemize}
\item{$K(s,n,p,D)$ is $n$-equivalent to $K$,}
\item{$\hat{Z}^{LMO}(\Sigma^p_{K(s,n,p,D)})-\hat{Z}^{LMO}(\Sigma^p_K) = p\alpha(s,n,p) D\ +\ $ terms of higher order.} 
\end{itemize}
\end{theorem}

\

See Example \ref{examp} 
for an illustration of the simple idea behind this theorem.

\begin{remark}
In some sense, leaving the $n$-triviality of the realising knot a free 
parameter is the hard part of the above theorem. It is likely, for example,
that such a realisation theorem (with $n$ fixed at $s$) is an extension
of Garoufalidis's investigations of the LMO invariant on branched cyclic 
covers over doubled knots. 

Note that if the $n$-triviality is not required to be free above, 
and $p>2$, then a realisation with $\alpha=1$ exists. Can we find an
improvement of the above theorem in the stated generality where $\alpha(s,n,p)=1$?
\end{remark}

A corollary of Theorem \ref{notft} is:

\begin{corollary}
Take a positive integer $p$.
If $v$ is a rational valued three-manifold invariant 
factoring non-trivially through the LMO invariant on integral homology 
spheres, then $v\circ \Sigma^p$ is not a 
finite-type invariant of knots.
\end{corollary} 

The proof of this corollary is indicated in Section \ref{corolproof}.

\subsection{Graph Y-links and the LMO invariant}

Graph Y-links in rational homology three-spheres realise particular first
non-vanishing terms of the LMO invariant. It is convenient here to cite
a direct calculation \cite{LeGr}. Below, let $\hat{Z}^{LMO}$ denote the 
version of the LMO invariant normalised by powers of the order of the first 
homology, following the notation of \cite{LeGr}.
\begin{theorem}
Let $M$ be a rational homology three-sphere and let $L$ be a 
connected graph Y-link in $M$ with some associated uni-trivalent diagram $D_L$ 
(i.e. represented by the underlying graph of $L$ with some choice of 
orientation at each vertex).
Then, 
\begin{equation}
\hat{Z}^{LMO}(M^L) - \hat{Z}^{LMO}(M) = \pm D_L + \mbox{terms of higher order}.
\end{equation}
\end{theorem}
Using Lemma \ref{minusintro}, one can always realise the opposite sign
by introducing a half-twist. 

\begin{example}\label{examp}

Before we consider the details of the realisation theorem, 
let us consider an example
to highlight the simple idea involved. Consider the following decoration
of the unknot $U$ by a graph Y-link $\kappa_n$:

\

\vspace{1.8cm}
\[
\thetaA
\]\vspace{1.5cm}

According to Theorem \ref{mainthb}, $\Sigma^p_{U^{\kappa_n}}-\Sigma^p_U$
lies in $\CF^Y_4\CM$, and moreover is proportional to the difference:

\

\[
\Picture{
\put(0,-4.75){
$\thetab$}}
\hspace{2cm} - \hspace{1cm} S^3\ \ \in \CG^Y_4\CM
\]
\vspace{0.5cm}

To calculate the proportionality constant, we follow the procedure described
in Example \ref{compexamp}. 
Start by introducing a variable for each leg $(\epsilon_1,\,\ldots,
\epsilon_n)$. If $\epsilon_j$ is zero, then that leg is to be sawn off; if
it is 1 then that leg is sawn off with an extra wrap around the knot (see
Example \ref{compexamp}).

The underlying trivalent graph has three cycles. Two of the cycles are 
unlinked from the knot for every $n$-tuplet $(\epsilon_1,\,\ldots,\epsilon_n)$.
The linking number of the third cycle with the knot will be 
$\epsilon_1+\ldots +\epsilon_n$.

Following Example \ref{compexamp}, letting
\[
F(n,t) = \sum_{(\epsilon_1,\,\ldots,\epsilon_n)=(0,\,\ldots,0)}^{(1,\,\ldots,1)}(-1)^{\epsilon_1+\ldots+\epsilon_n}
t^{\epsilon_1+\ldots+\epsilon_n} = (1-t)^n,
\]
then
\[
\Sigma^p_{U^{\kappa_n}} - S^3 = \left(
\sum_{q=0}^{p-1}F\left(n,e^{\frac{2\pi i q}{p}}\right)
\right) 
((S^3)^{\lambda_{\kappa_n}} - S^3 )\ \in \CG^Y_4\CM.
\]
Thus
\[
\hat{Z}^{LMO}(\Sigma^p_{U^{\kappa_n}}) = 1 + \epsilon(\sum_{q=0}^{p-1}( 1 - e^{\frac{2\pi i q}{p}})^n) \kappa + \mbox{terms of higher degree},
\]
for some $\epsilon = \pm 1$.

\end{example}

Given a knot $K$, an oriented trivalent diagram $D$, and a positive
integer $n$, 
we {\it choose} a knot $K_{(D,n)}$ as follows. 
Select a graph Y-link $L$ in the 
three-sphere so that 
\[
\hat{Z}^{LMO}( (S^3)^{L} - S^3 )= D\ +\ \mbox{terms of higher degree}.
\]
Then, locate $L$ in a ball disjoint from $K$ in the three-sphere.
Select some edge of $L$, give it
an orientation, and add $n$ legs joining that edge to the knot $K$ in 
that ball, so that each edge is added according to the following 
orientation (but otherwise, in any fashion):

\

\vspace{0.25cm}

\hspace{6cm}\addlegs\hspace{2.25cm}
\vspace{2cm}

There are many choices in this definition. Nevertheless:
\begin{lemma}
\[
\hat{Z}^{LMO}( \Sigma^p_{K_{(D,n)}} - \Sigma^p_K )= 
\left(
\sum_{q=0}^{p-1}F\left(n,e^{\frac{2\pi i q}{p}}\right)
\right) 
D\ +\ \mbox{terms of higher degree}.
\]
\end{lemma}
This follows from the same logic as in the above example. It is clear
that if $D$ has surplus $s$, then $K_{(D,n)}$ is obtained from
$K$ by surgery on a graph Y-link of degree $\frac{s}{2}+n$, and so
is $(\frac{s}{2}+n)$-equivalent to $K$ \cite{Hab}.

\subsection{Proof of Theorem \ref{notft}}

Without loss of generality, let the combination to be realised 
be written $S=\sum_{i=1}^{t}D_i$. Take some positive integer $l$.
Let $K_{(S,l)}$ denote a knot:
\[
(\ldots((K_{(D_1,l)})_{(D_2,l)})\ldots )_{(D_t,l)}.
\]
It is clear that
\[
\hat{Z}^{LMO}( \Sigma^p( K_{(S,l)} - K ) )
=
\left(
\sum_{q=0}^{p-1}F\left(l,e^{\frac{2\pi i q}{p}}\right)
\right) S + \mbox{terms of higher order}. 
\]

The number of legs to be added, $l$, is a free parameter, so that
one can make the knot $K_{(S,l)}$ $n$-equivalent to $K$ for arbitrarily
large $n$.

The remaining difficulty is covered by the following lemma.

\rtb

\begin{lemma}
\[
\mbox{lim}_{l\rightarrow \infty} 
\left(
\sum_{q=0}^{p-1}F\left(l,e^{\frac{2\pi i q}{p}}\right)
\right) \neq 0.
\]
\end{lemma}

\underline{Proof}

Let $\omega = e^{\frac{2\pi i}{p}}.$
Define a function $f(l)= 
\sum_{q=0}^{p-1}F(l,\omega^q)$. 
Take some positive integer $l$. We will show that $f(l')$ is
non-zero for at least 
one $l'$ such that $l\leq l' < l+p$.

We assume that $f(l+i)=0$ for $0\leq i < p$, and derive a 
contradiction. Now, we have assumed that $f(l)=0$, so that,
letting $T_q = (1-\omega^{q})$,
 
\[
f(l+1)  =  \sum_{q=0}^{p-1}(1-\omega^{q})T_q^{l} 
= -\sum_{q=0}^{p-1}\omega^{q}T_q^{l} = 0.
\]
Proceeding with this assumption, one finds that
\[
f(l+j) 
 = 
(-1)^{j}\sum_{q=0}^{p-1}\omega^{jq}T_q^{l} = 0.
\]
In other words, under this assumption, we have a matrix equation:
\[
\left(
\begin{array}{r}
f(l) \\ -f(l+1) \\ f(l+2) \\ . \\ . \\ . \\ (-1)^{p-1}f(l+p-1)
\end{array}
\right) =
\left(
\begin{array}{lllllll}
1 & 1 & 1 & . & . & . & 1 \\
1 & \omega & \omega^2 & . & . & . & \omega^{p-1} \\
1 & \omega^2 & \omega^4 & . & . & . & \omega^{2(p-1)} \\
. & . & . & . & . & . & . \\
. & . & . & . & . & . & . \\
. & . & . & . & . & . & . \\
1 & \omega^{(p-1)} & \omega^{2(p-1)} & . & . & . & \omega^{(p-1)(p-1)}
\end{array}
\right)
\left(
\begin{array}{l}
T_0^l \\ T_1^l \\ T_2^l \\ . \\ . \\ . \\ T_{p-1}^l
\end{array}
\right)\ \ \ \ \ \ \ \ \ \ \ \ 
\]
The vector on the right in clearly non-zero, and the matrix is a Vandermonde
matrix, and is non-singular, so the vector on the left cannot be zero.
This is our contradiction. 
\rtb 

\

\subsection{Proof of Corollary \ref{thecorol} }\label{corolproof}

\

This is equivalent to finding for {\it any} vector $D\in \CA_s$ 
an element $E_{(D,n)} \in \CF_n\CK$ such that
\[
\KI^{LMO}(\Sigma^p(E_{(D,n)})) = D + \mbox{higher order terms}.
\]

This follows easily from the Theorem \ref{notft} (freely multiplying by 
rationals) and the observation that for knots $K$ and $L$:
\begin{equation}\label{multeqn}
\Sigma^p( K\# L) = \Sigma^p(K) \# \Sigma^p(L).
\end{equation}

\rtb

\

\appendix
\section{The Goussarov-Habiro theory of finite-type invariants of three-manifolds}\label{GHapp}\label{GHT}

\

The theory to be recalled in this section is due, independently, to 
Goussarov \cite{G3}
and to Habiro \cite{Hab}. The original finite type theory (in the setting
of $\Zset HS^3$s) is due to 
Ohtsuki \cite{Oht}; in that setting there are also definitions due to 
Garoufalidis and Levine \cite{GL}; and the first definition to consider
the set of all three-manifolds was due to Cochran and 
Melvin \cite{CM}. 
The Goussarov-Habiro theory is defined in terms of 
a move on three-manifolds due to Matveev \cite{Mat}, closely related
to a move of Murakami and Nakanishi's, on links \cite{MN}. 
Our summary will be concise: we expect
expositions of this theory 
to appear in the future (\cite{GGP}).

A {\it clasper} is an embedding of a standard 
band-summed pair of annuli into a three-manifold. 
The annuli are called the {\it leaves} of the clasper.
This gives placement information for a two-component
framed link, as follows, and a {\it move} on that clasper replaces that
manifold with the manifold obtained by doing surgery 
on that link. The tubes in the diagram below can contain
part of the surgery link presenting the three-manifold, and any
other objects that may be embedded in the three-manifold: \vspace{0.1cm}

\[
\setlength{\unitlength}{30pt}
\hspace{1.3cm}
\Picture{
\qbezier(1.1,-0.5)(1.1,-0.8)(1.1,-1.4)
\qbezier(0.9,1)(0.9,0.5)(0.9,-0.3)
\qbezier(0.9,-0.5)(0.9,-0.8)(0.9,-1.4)
\qbezier(1.1,1)(1.1,0.5)(1.1,-0.3)
\qbezier(-0.9,1)(-0.9,0.5)(-0.9,0.1)
\qbezier(-0.9,-1.4)(-0.9,-1.0)(-0.9,-0.1)
\qbezier(-1.1,1)(-1.1,0.5)(-1.1,0.1)
\qbezier(-1.1,-1.4)(-1.1,-1.0)(-1.1,-0.1)
\qbezier(1.1,0.1)(1.4,0.1)(1.4,-0.2)
\qbezier(1.1,-0.5)(1.4,-0.5)(1.4,-0.2)
\qbezier(1.1,-0.1)(1.2,-0.1)(1.2,-0.2)
\qbezier(1.1,-0.3)(1.2,-0.3)(1.2,-0.2)
\qbezier(1.1,-0.5)(0.9,-0.5)(0.9,-0.5)
\qbezier(1.1,-0.3)(0.9,-0.3)(0.9,-0.3)
\qbezier(0.9,-0.5)(0.6,-0.5)(0.6,-0.2)
\qbezier(0.9,-0.3)(0.8,-0.3)(0.8,-0.2)
\qbezier(0.9,-0.1)(0.8,-0.1)(0.8,-0.2)
\qbezier(0.9,0.1)(0.6,0.1)(0.6,-0.2)
\qbezier(-1.1,-0.5)(-1.4,-0.5)(-1.4,-0.2)
\qbezier(-1.1,0.1)(-1.4,0.1)(-1.4,-0.2)
\qbezier(-1.1,-0.3)(-1.2,-0.3)(-1.2,-0.2)
\qbezier(-1.1,-0.1)(-1.2,-0.1)(-1.2,-0.2)
\qbezier(-1.1,0.1)(-0.9,0.1)(-0.9,0.1)
\qbezier(-1.1,-0.1)(-0.9,-0.1)(-0.9,-0.1)
\qbezier(-0.9,0.1)(-0.6,0.1)(-0.6,-0.2)
\qbezier(-0.9,-0.1)(-0.8,-0.1)(-0.8,-0.2)
\qbezier(-0.9,-0.3)(-0.8,-0.3)(-0.8,-0.2)
\qbezier(-0.9,-0.5)(-0.6,-0.5)(-0.6,-0.2)
\qbezier(0.6,-0.1)(-0,-0.1)(-0.6,-0.1)
\qbezier(0.6,-0.3)(0,-0.3)(-0.6,-0.3)
}
\ \ \ \ \ \ \ \ \rightsquigarrow \hspace{3cm} 
\Picture{
\qbezier(1.1,-0.5)(1.1,-0.8)(1.1,-1.4)
\qbezier(0.9,1)(0.9,0.5)(0.9,-0.3)
\qbezier(0.9,-0.5)(0.9,-0.8)(0.9,-1.4)
\qbezier(1.1,1)(1.1,0.5)(1.1,-0.3)
\qbezier(-0.9,1)(-0.9,0.5)(-0.9,0.1)
\qbezier(-0.9,-1.4)(-0.9,-1.0)(-0.9,-0.1)
\qbezier(-1.1,1)(-1.1,0.5)(-1.1,0.1)
\qbezier(-1.1,-1.4)(-1.1,-1.0)(-1.1,-0.1)
\qbezier(-1.1,0)(-1,0)(-0.1,0)
\qbezier(-0.1,0)(0.1,0)(0.1,-0.2)
\qbezier(1.1,-0.4)(0.1,-0.4)(0.1,-0.4)
\qbezier(0.1,-0.4)(-0.1,-0.4)(-0.1,-0.2)
\qbezier(0.1,-0.2)(0.1,-0.275)(0.05,-0.35)
\qbezier(-0.1,-0.2)(-0.1,-0.125)(-0.05,-0.05)
\qbezier(0.1,0)(0.1,0)(0.9,0)
\qbezier(-0.1,-0.4)(-0.1,-0.4)(-0.9,-0.4)
\qbezier(1.1,0)(1.3,0)(1.3,-0.2)
\qbezier(1.1,-0.4)(1.3,-0.4)(1.3,-0.2)
\qbezier(-1.1,0)(-1.3,0)(-1.3,-0.2)
\qbezier(-1.3,-0.2)(-1.3,-0.4)(-1.1,-0.4)
\put(0.4,0.15){\mbox{0}}
\put(-0.6,0.15){\mbox{0}}
}
\]
\vspace{1.4cm}

In this work, claspers in $S^3$ will be depicted via blackboard 
framed diagrams. That is, circles should be thickened to annuli, and
the other arcs should be thickened to bands, in the plane of 
the diagram. The following
symbols will be used to indicate half-twists of bands:\vspace{0.1cm}

\[
\Picture{
\put(0,0){\circle{1}}
\qbezier(0,0.5)(0,2)(0,2)
\qbezier(0,-0.5)(0,-2)(0,-2)
\qbezier(0.35,0.35)(0,0)(-0.35,-0.35)
}\ \ \mbox{\em represents}\ \ \ \ \ \ \ \ 
\Picture{
\qbezier(0.4,2)(0.4,1.2)(0.4,0.8)
\qbezier(0.4,0.8)(0.4,0.4)(0,0)
\qbezier(-0.4,-0.8)(-0.4,-0.4)(0,0)
\qbezier(-0.4,-2)(-0.4,-1.2)(-0.4,-0.8)
\qbezier(-0.4,2)(-0.4,0.4)(-0.2,0.2)
\qbezier(0.4,-2)(0.4,-0.4)(0.2,-0.2)
} \ \ \mbox{\em and}\ \ \ \ \ \  
\Picture{
\put(0,0){\circle{1}}
\qbezier(0,0.5)(0,2)(0,2)
\qbezier(0,-0.5)(0,-2)(0,-2)
\qbezier(-0.35,0.35)(0,0)(0.35,-0.35)
}\ \ \mbox{\em represents}\ \ \ \ \ \ 
\Picture{
\qbezier(-0.4,2)(-0.4,1.2)(-0.4,0.8)
\qbezier(-0.4,0.8)(-0.4,0.4)(0,0)
\qbezier(0.4,-0.8)(0.4,-0.4)(0,0)
\qbezier(0.4,-2)(0.4,-1.2)(0.4,-0.8)
\qbezier(0.4,2)(0.4,0.4)(0.2,0.2)
\qbezier(-0.4,-2)(-0.4,-0.4)(-0.2,-0.2)
}.
\]\vspace{0.5cm}

A {\it Y-component} is an embedding of a standard triple of annuli band-summed
into a disk, into a three-manifold. This gives placement information
for a six-component surgery link. A {\it move} on that Y-component
replaces that manifold with the manifold obtained by doing surgery 
on that link. This is Matveev's Borromeo move and Murakami-Nakanishi's
Delta-unknotting move. 
It is convenient, here, to build a standard Y-component
from three claspers, as follows:\vspace{2cm}

\[
\hspace{0.8cm}\ycomp\hspace{2cm}\rightsquigarrow\hspace{3cm}\actycomp
\]
\vspace{2cm}

\begin{flushleft}
where the following convention has been used:
\end{flushleft}

\[
\epsfxsize = 2.5cm
\fig{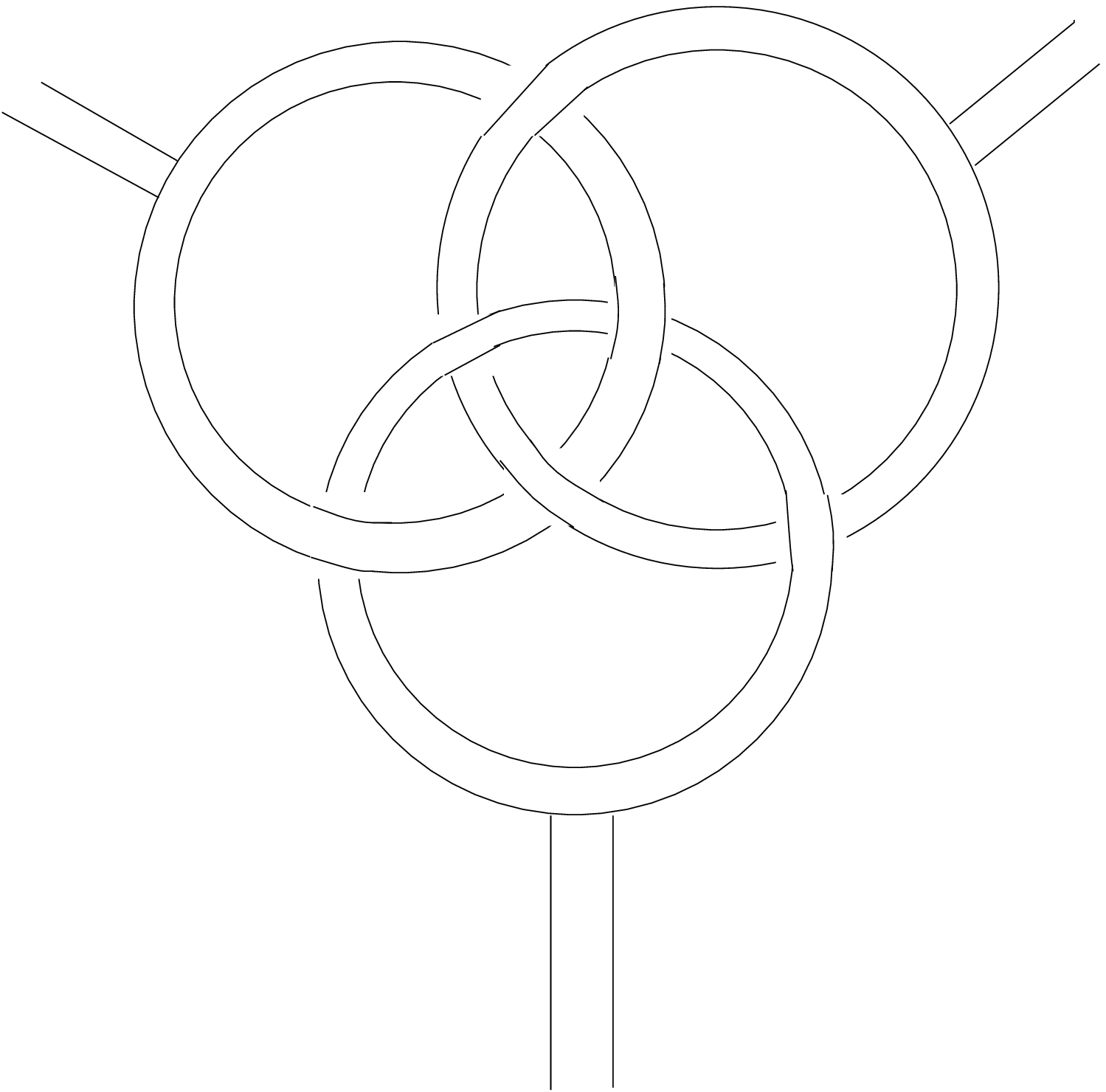}\ \ \ \mbox{\em is indicated by}\ \ \ \ \ \ \ \ \
\Picture{
\put(0.1,0){\circle*{1}}
\qbezier(-0.14,0.14)(-2,2)(-2,2)
\qbezier(0.14,0.14)(2,2)(2,2)
\qbezier(0,-0.14)(0,-2)(0,-2)
}\ \ \ \ .
\]

It can be shown that this
move generates an equivalence relation on the set of closed
oriented three-manfolds. 

\begin{theorem}[\cite{Mat}]\label{Mattheo}
For two closed oriented three-manifolds $M$ and $N$,
there is an isomorphism $H_1(M) \rightarrow H_1(N)$ preserving the
linking form on the torsion if and only if $N$ is obtained from $M$
via a finite sequence of Y-moves.
\end{theorem}

The Goussarov-Habiro
theory is based on a filtration of $\CM$, 
the free abelian group generated by homeomophism classes of 
closed oriented three-manifolds. 

The definition uses a {\it Y-link} in a three-manifold, which is 
a collection of disjointly embedded Y-components.
Take a $\mu$-component Y-link $L$, and let
$M_L$ denote the three-manifold obtained by performing the moves on
the Y-components of $L$.
Define a vector $[M,L] \in \CM$ corresponding to a pair of some $M$ and some
Y-link $L$ in $M$ by
\[
[M,L] = \sum_{L'\subset L}(-1)^{\#L -\# L'} M_{L'},
\]
where the sum is over all Y-sublinks of $L$ (so that there will be
$2^\mu$ such). 

\begin{definition}\label{elemdef}\label{GHD}
Define the subspace $\CF_n^Y \CM \subset \CM$ to be the subspace spanned
by all vectors $[M,L]$ with $L$ an $n$-component Y-link.
\end{definition}

It is clear that this filters $\CM$,
\[
\CM \supset \CF^Y_1\CM
\supset \CF^Y_2\CM
\supset \CF^Y_3\CM
\ldots\ ,
\]
and that a theory of finite-type invariants of three-manifolds can be 
constructed by declaring an invariant to be finite-type of order $\geq n$ if it 
vanishes on the subspace $\CF^Y_{n+1}\CM$. 

The associated graded quotients $\frac{\CF^Y_{n}\CM}{\CF^Y_{n+1}\CM}$,
denoting them $\CG^Y_{n}\CM$, are finite dimensional, with
a finite spanning set of generators of the form $[M,L]$ where $L$
is a {\it graph Y-link}.
These are Y-links where
every leaf of every Y-component links precisely one other leaf
in a ball, as follows:\vspace{1.5cm}

\hspace{5.5cm}\seifcrossb\vspace{1.8cm}

\hspace{-0.5cm}
or is some standard link fixed to represent some element of homology. The
full class will not occur in this work, so when we refer to a 
graph Y-link, we will mean specifically a Y-link where every leaf meets
another in a ball.

To a graph Y-link is associated a uni-trivalent graph by 
taking a trivalent vertex for every Y-component, and joining edges
when the associated leaves link in a ball. If the graph associated
to some graph Y-link is connected, call the associated graph Y-link
{\it connected}. 
\begin{lemma}\label{subgr}
If $L$ is a connected graph Y-link, and $L'$ is a Y-sublink
of $L$ not $L$, then $M_{L'} \simeq M$. A consequence is that
\[
[M,L'] = 0.
\]
\end{lemma}
This follows from the fact that any clasper or Y-component with a
leaf that bounds a disc in the complement of the rest of the Y-link,
may be removed without affecting the result:

\[
\setlength{\unitlength}{20pt}
\Picture{
\qbezier(-1.1,1)(-1.1,0.5)(-1.1,0.1)
\qbezier(-1.1,-1.4)(-1.1,-1.0)(-1.1,-0.1)
\qbezier(-0.9,1)(-0.9,0.5)(-0.9,0.1)
\qbezier(-0.9,-1.4)(-0.9,-1.0)(-0.9,-0.1)
\qbezier(1.1,0.1)(1.4,0.1)(1.4,-0.2)
\qbezier(1.1,-0.5)(1.4,-0.5)(1.4,-0.2)
\qbezier(1.1,-0.1)(1.2,-0.1)(1.2,-0.2)
\qbezier(1.1,-0.3)(1.2,-0.3)(1.2,-0.2)
\qbezier(1.1,-0.5)(0.9,-0.5)(0.9,-0.5)
\qbezier(1.1,-0.3)(0.9,-0.3)(0.9,-0.3)
\qbezier(0.9,-0.5)(0.6,-0.5)(0.6,-0.2)
\qbezier(0.9,-0.3)(0.8,-0.3)(0.8,-0.2)
\qbezier(0.9,-0.1)(0.8,-0.1)(0.8,-0.2)
\qbezier(0.9,0.1)(0.6,0.1)(0.6,-0.2)
\qbezier(-1.1,-0.5)(-1.4,-0.5)(-1.4,-0.2)
\qbezier(-1.1,0.1)(-1.4,0.1)(-1.4,-0.2)
\qbezier(-1.1,-0.3)(-1.2,-0.3)(-1.2,-0.2)
\qbezier(-1.1,-0.1)(-1.2,-0.1)(-1.2,-0.2)
\qbezier(-1.1,0.1)(-0.9,0.1)(-0.9,0.1)
\qbezier(0.9,0.1)(0.9,0.1)(1.1,0.1)
\qbezier(1.1,-0.1)(0.9,-0.1)(0.9,-0.1)
\qbezier(-1.1,-0.1)(-0.9,-0.1)(-0.9,-0.1)
\qbezier(-0.9,0.1)(-0.6,0.1)(-0.6,-0.2)
\qbezier(-0.9,-0.1)(-0.8,-0.1)(-0.8,-0.2)
\qbezier(-0.9,-0.3)(-0.8,-0.3)(-0.8,-0.2)
\qbezier(-0.9,-0.5)(-0.6,-0.5)(-0.6,-0.2)
\qbezier(0.6,-0.1)(-0,-0.1)(-0.6,-0.1)
\qbezier(0.6,-0.3)(0,-0.3)(-0.6,-0.3)
}
\ \ \sim \ \ \ \ \ \ \ \ \ \  
\Picture{
\qbezier(-1.1,1)(-1.1,0.5)(-1.1,-1.4)
\qbezier(-1.1,-1.4)(-1.1,-1.0)(-1.1,-1.4)
\qbezier(-0.9,1)(-0.9,0.5)(-0.9,-1.4)
\qbezier(-0.9,-1.4)(-0.9,-1.0)(-0.9,-0.1)
}\hspace{-1.4cm}.
\]\vspace{0.7cm}

In the remainder of this section, we describe some aspects of
the manipulation of clasper graphs that will be useful in
ensuing calculations.

In this work, we will encounter situations where it is convenient
to present a Y-link by decorating some other Y-link with a collection of
claspers: the desired Y-link being recovered by surgery on those claspers.
So, on occasions when precision is required, 
we will call such a link a {\it mixed Y-link}. The
definition of the vector $[M,L]$ is extended to mixed Y-links,
taking the alternating sum over the Y-components only.
A {\it sublink} of a mixed link, is obtained by forgetting a number
of Y-components or claspers.

Now we collect some moves that will be useful in this paper.

\begin{lemma}\label{insertclasp}
The result of surgery on two mixed Y-links which differ in a ball as follows
is the same.
\ \vspace{0.5cm}\

\hspace{2.4cm}\brek\hspace{4cm}\brekb
\end{lemma}\vspace{0.25cm}

\begin{lemma}\label{legbreak}
The result of surgery on two mixed Y-links which differ in a ball as follows
is the same.
\ \vspace{1.5cm}\

\hspace{2cm}\legccp\hspace{5.25cm}\Picture{\put(0,-4){$\leggedy$}}
\end{lemma}\vspace{3.2cm}
\begin{lemma}\label{add}
Let $L_A$, $L_B$ and $L_C$ be $n$-component 
Y-links that differ in a ball as follows. The dashed
part of the leaf indicates that that part follows some path in the 
three-manifold before returning to the ball in question.
\vspace{2.5cm}

\hspace{3.6cm}
\reltripa\hspace{3.2cm}\reltripb\hspace{3.2cm}\reltripc\hspace{3cm}
\vspace{0.1cm}

\[
[M,L_A] = [M,L_B] + [M,L_C] \in \CG^Y_n\CM.
\]
\end{lemma}

\newpage

\begin{lemma}
\label{crosstube}

Let the Y-links $L_A$ and $L_B$ differ as follows. The tube can contain 
parts of surgery components or part of the rest of the graph:
\vspace{2cm}

\hspace{3cm}\legppc\hspace{4.5cm}\legpppc
\vspace{3cm}

\[
[M,L_A] = [M,L_B] \in \CG^Y_n\CM.
\]
\end{lemma}

Finally, let $L_A$ and $L_B$ differ in a ball as follows:\vspace{1cm}

\hspace{3cm}\legc\hspace{4.5cm}\legcc
\vspace{3.5cm}

\begin{lemma}
\begin{equation}\label{minusintro}
[M,L_A] = -[M,L_B] \in \CG^Y_n\CM.
\end{equation}
\end{lemma}


\subsection{n-equivalence}\label{nequivdefn}

The property of $n$-equivalence was
introduced by Ohyama in the setting of knots \cite{Ohy};
its importance for Vassiliev theory
was observed by Goussarov \cite{G1}.

If two three-manifolds are $n$-equivalent then
their difference lies in $\CF^Y_{n+1}\CM$ and on this pair
all finite-type invariants of order less than or equal to $n$ agree.

\begin{definition}
A {\it $n+1$-scheme} for $M$ is a mixed Y-link $L$ in $M$ together with 
a set of $n+1$ disjoint Y-sublinks of $L$, $L_1$ up to $L_{n+1}$. 
\end{definition}
For some $n+1$-tuplet $\{i_1,\ldots,i_{n+1}\}$, 
where $i_k$ is either 0 or 1, the
notation $L_{i_1,\ldots,i_j}$ is used to denote the Y-link
that is obtained by forgetting those sublinks whose associated index is $1$.

\begin{definition}
A three-manifold $N$ is {\it $n$-equivalent} to a three-manifold $M$
if $M$ has an $n+1$-scheme $\{L; L_1,\ldots L_{n+1}\}$ such that 
\begin{itemize}
\item{$M_{L_{0,\ldots,0}}\simeq N$,}
\item{$M_{L_{i_1,\ldots,i_{n+1}}}\simeq M$ for any other multiplet.}
\end{itemize}
\end{definition}

In such a situation we will say that $\{L; L_1,\ldots L_{n+1}\}$ is
an $n+1$-scheme {\it relating} $N$ to $M$.

\begin{lemma}\label{nequivlem}
If $M$ is $n$-equivalent to $N$, then
\[
M - N \in \CF^Y_{n+1}\CM.
\]
\end{lemma}
 
In such a situation we also have a nice expression for $M - N$ in the
graded space $\frac{\CF^Y_{n+1}\CM}{\CF^Y_{n+2}\CM}$. 

To introduce this expression we need to be more specific with some notation.
Denote the relating $n+1$-scheme in $M$ by $\{L ; L_1,\ldots,L_{n+1}\}$. 
Let $\sigma(i)$ be the function 
giving the number of Y-components of the Y-sublink $L_i$.
Order the Y-components of each Y-sublink $L_i$.
For an $n+1$-tuple $(a_1,\ldots,a_{n+1})$, where $1\leq a_i\leq \sigma(i)$, let
$L^{(a_1,\ldots,a_{n+1})}$ be the mixed Y-link 
in $M$ obtained by forgetting all Y-components of each of these
Y-sublinks, except precisely one Y-component from each $L_i$: that
is, from $L_i$ choose $a_i$.

\begin{lemma}\label{niceexprlem}
\[
N - M = \sum_{(a_1,\ldots,a_{n+1})=(1,\ldots,1)}^{(\sigma(1),\ldots,\sigma(n+1))}
[M,L^{(a_1,\ldots,a_{n+1})}] \in \CG^Y_{n+1}\CM.
\]
\end{lemma}

\bibliographystyle{amsalpha}

\begin{thebibliography}{MS2}
\bibitem[CM]{CM} T. Cochran and P. Melvin, {\em Finite type invariants
of 3-Manifolds}, to appear in Inventiones.
\bibitem[Dav1]{Dav1} A. Davidow, {\em Casson's invariant and iterated torus
knots}, Knots '90, Walter de Gruyter, 151--161.
\bibitem[Dav2]{Dav2} A. Davidow, {\em Casson's invariant and twisted double knots},
Top. Appl. {\bf 58} (1994) 93--101.
\bibitem[F]{F} R.H. Fox, {\em Free differential calculus III}, Ann. of Math. {\bf 64} (1956), 407--419.
\bibitem[Gar]{Gar} S. Garoufalidis, 
{\em Signatures of links and finite-type invariants of cyclic branched
covers}, preprint 1999.
\bibitem[GH]{GH} S. Garoufalidis and N. Habegger, {\em 
The Alexander polynomial
and finite type 3-manifold invariants}, to appear in Math. Annalen.
\bibitem[GGP]{GGP} S. Garoufalidis, M. Goussarov and M. Polyak, 
{\em Topological calculus of $y$-graphs and equivalence of finite
type invariants}, in preparation.
\bibitem[GL]{GL} S. Garoufalidis and J. Levine, 
{\em Finite type 3-manifold invariants and the structure of the 
Torelli group I}, to appear in the Journal of Diff. Geom.
\bibitem[GL2]{GL2} S. Garoufalidis and J. Levine, 
{\em Tree level invariants
of three-manifolds, Massey products and the Johnson homomorphism},
preprint 1999.
\bibitem[G]{G} L. Goeritz, {\em Knotten und quadratische Formen}, Math. Zeit.
{\bf 36} (1933), 647--654.
\bibitem[Gor]{Gor} C.McA. Gordon, {\em Knots whose branched cyclic coverings
have periodic homology}, Trans. Amer. Math. Soc. {\bf 168} (1972), 357--370. 
\bibitem[G1]{G1} M. Goussarov, {\em N-equivalent knots and invariants of
finite degree}, Topology of manifolds and varieties, ed: Viro, O., Amer. Math. Soc. (1994) 173--192.
\bibitem[G2]{G2} M. Goussarov, {\em Knotted graphs and a geometrical 
technique of $n$-equivalence}, POMI Sankt Petersburg preprint, circa 1995. 
\bibitem[G3]{G3} M. Goussarov, {\em Finite-type invariants and
$n$-equivalence of 3-manifolds}, to appear in Comptes Rendus.
\bibitem[GW]{GW} Papers of M. Goussarov, {\tt http://www.ma.huji.ac.il/~drorbn/Goussarov}, maintained by D. Bar-Natan.
\bibitem[HB]{HB} N. Habegger and A. Beliakova, {\em The Casson-Walker-Lescop
Invariant as a Quantum 3-manifold Invariant}, preprint 1997.
\bibitem[Hab]{Hab} K. Habiro, {\em Clasper theory and finite-type 
invariants of links}, Geometry and Topology, {\bf 4}, (2000), 1--83.
\bibitem[HT]{HT} K. Habiro, {\em Aru musubime no kyokusyo so usa no zoku ni tsuite}, Tokyo University Master's Thesis (1994).
\bibitem[HK]{HK} F. Hosokawa and S. Kinoshita, {\em On the homology
group of branched cyclic covering spaces of links}, Osaka Math. J., {\bf 12},
(1960), 331--335.
\bibitem[Hos]{Hos} J. Hoste, {\em The first coefficient of the Conway polynomial},
Proc. Amer. Math. Soc. {\bf 95} (1985) 299--302.
\bibitem[Ish]{Ish} K. Ishibe, {\em The Casson-Walker invariant for
branched cyclic covers of $S^3$ branched over a doubled knot},
Osaka J. Math. {\bf 34} (1997) 481--495
\bibitem[Kan]{Kan} T. Kanenobu, {\em Examples on polynomial invariants
of knots and links}, Math. Ann., {\bf 275}, (1986), 555--572
\bibitem[K]{K} A. Kricker, {\em Covering spaces over claspered knots},
preprint 1998 (this has appeared in the informal proceedings of the workshop
``Knot Theory'', Kyoto, October 1998, ed. M. Sakuma).
\bibitem[KGr]{KGr} A. Kricker, {\em Clasper moves and the rows of the coloured
Jones function}, seminar at the Ecole d'Ete, University Joseph Fourier, Grenoble, July 1999.
\bibitem[LMMO]{LMMO} 
T. Le, J. Murakami, H. Murakami and T. Ohtsuki,
{\em A three-manifold invariant via the Kontsevich integral},
to appear in Osaka Journal of Math.
\bibitem[LeGr]{LeGr} T.Q.T. Le, {\em The LMO invariant}, 
Notes accompanying lectures at the Ecole d'Ete, University Joseph Fourier, Grenoble, June, 1999.
\bibitem[LMO]{LMO} T.T.Q. Le, J. Murakami. and Ohtsuki., T.,
{\em A universal quantum invari1ant of three-manifolds}, Topology, {\bf 37}, (1998),
539--574.
\bibitem[Les]{Les} C. Lescop, {\em Global surgery formula for the Casson-Walker
invariant }, Annals of Math Studies No 140, Princeton University.
\bibitem[Lic]{Lic} R. Lickorish, {\em An introduction to knot theory},
Graduate Texts in Mathematics, {\bf 175}, (1997).
\bibitem[L]{L} J. Lieberum, {\em The LMO-invariant of 3-manifolds of
rank one and the Alexander polynomial}, preprint January 2000
\bibitem[Mat]{Mat} 
S.W. Matveev, 
{\em Generalised surgery of three-dimensional manifolds and representations of homology spheres}, Math. Notices. Acad. Sci. USSR, {\bf 42:2} (1987) 651--656.
\bibitem[MN]{MN}
H. Murakami and Y. Nakanishi,
{\em On a certain move generating link homology}, Math. Ann. {\bf 284}, (1989),
75--89. 
\bibitem[MO]{MO}
H. Murakami, and T. Ohtsuki,
{\em Finite-type invariants of knots via their 
Seifert matrices}, Tokyo Institute of Technology and Waseda University
preprint 1998.
\bibitem[Mul1]{Mul1}
D. Mullins, {\em The generalised Casson invariant for 2-fold branched covers
of $S^3$ and the Jones polynomial}, Topology, {\bf 32}, (1993), 419-438.
\bibitem[Mul2]{Mul2}
D. Mullins, {\em The Casson invariant for two-fold branched cyclic covers of links}, Quantum Topology, Series on Knots and Everything {\bf 3}, World Scientific, 221--229.
\bibitem[NS]{NS}
S. Naik and T. Stanford, 
{\em A move on diagrams that generates
S-equivalence of knots}, preprint, 1999. 
\bibitem[Ng]{Ng}
K.Y. Ng,
{\em Groups of ribbon knots},
Topology {\bf 37}, (1998), 441--458.
\bibitem[Oht]{Oht}
T. Ohtsuki,
{\em A filtration of the set of integral homology 3-spheres}, Proc. of the
Int. Congress of Math., Vol. II, Berlin, (1998), 473--482.
\bibitem[Ohy]{Ohy}
Y. Ohyama, 
{\em A new numerical invariant of knots induced from their regular diagrams},
Top. Appl. {\bf 37}, (1990), 249-255.
\bibitem[Rol]{Rol}
D. Rolfsen,
{\em Knots and links}, Mathematics Lecture Series {\bf 7}, (1976), Publish
or Perish, inc.
\bibitem[Roz]{Roz}
L. Rozansky,
{\em A rational structure of generating functions for Vassiliev invariants},
Notes accompanying lectures at the Ecole d'Ete, University Joseph Fourier, Grenoble,
June 1999.
\bibitem[S]{S} T. Stanford, {\em Vassiliev invariants and knots modulo pure
braid subgroups}, United States Naval Academy preprint, 1998.
\bibitem[ST]{ST}
T. Stanford and R. Trapp,
{\em On knot invariants which are not of finite type},
preprint, March 1999.
\end{thebibliography}

\end{document}